\documentclass[11pt]{article}
\usepackage{amsthm,amssymb,amsmath}
\usepackage{graphicx,algorithmic,algorithm}%
\usepackage{enumerate}
\usepackage{tikz,tikz-cd}
\usepackage{tabularx}
\usepackage{booktabs}
\usepackage{multirow}
\usepackage{makecell}
\usepackage{csquotes}

\tikzcdset{scale cd/.style={every label/.append style={scale=#1},
    cells={nodes={scale=#1}}}}
\usetikzlibrary{shapes,fit}
\usetikzlibrary{arrows.meta}
\usepackage[letterpaper, margin=1in]{geometry}
\usepackage{hyperref}
\usepackage[capitalize]{cleveref}
\usepackage{todonotes}
\usepackage[normalem]{ulem}
\hypersetup{
	pdftitle=   {},
	pdfauthor=  {}
}
\usepackage{mathabx}
\usepackage{caption}
\usepackage{subcaption}
\usepackage{commath}

\crefname{figure}{Figure}{Figures}
\crefname{equation}{Equation}{Equations}

\newtheorem{theorem}{Theorem}[section]
\newtheorem{lemma}[theorem]{Lemma}
\newtheorem{observation}[theorem]{Observation}

\newtheorem{corollary}[theorem]{Corollary}

\newtheorem{claim}[theorem]{Claim}

\newtheorem*{thm:main1}{Theorem~\ref{thm:main1}}
\newtheorem*{thm:main2}{Theorem~\ref{thm:main2}}
\newenvironment{clproof}{\begin{list}{}{%
			\setlength{\leftmargin}{3mm}%
		} \item {\it Proof.} }{\hfill$\lozenge$\end{list}}

\newcommand\extrafootertext[1]{%
    \bgroup
    \renewcommand\thefootnote{\fnsymbol{footnote}}%
    \renewcommand\thempfootnote{\fnsymbol{mpfootnote}}%
    \footnotetext[0]{#1}%
    \egroup
}

\DeclareMathOperator{\mim}{mim}

\DeclareMathOperator{\tww}{tww}

\newcommand{\ce}{\text{CE}}
\newcommand{\cc}{\text{CC}}
\newcommand{\cs}{\text{CS}}

\title{Graphs of twin-width at most~2 have bounded mim-width}
\author{%
  \'{E}douard Bonnet\thanks{CNRS, ENS de Lyon, Université Claude Bernard Lyon 1, LIP UMR 5668, Lyon, France.\newline
  The work of É.B. was supported by the ANR project OPAL (ANR-26-CE48-5754).}
  \and
  Yeonsu Chang\thanks{Department of Mathematics, Hanyang University, Seoul, South Korea.}%
}
\date{}

\begin{document}

\maketitle

\begin{abstract}
  We show that every graph of twin-width at most~$2$ has mim-width at most~$8$.
  This subsumes the facts that weakly sparse graphs of twin-width at~most~$2$ have bounded treewidth [Bergougnoux et al., ISAAC '23], and edge-stable graphs of twin-width at~most~$2$ have bounded rank-width [Heinrich et al., WG '26].  
Given a~\mbox{$2$-sequence}, we construct a~branch decomposition of mim-width~$8$ in polynomial time.
This gives polynomial-time algorithms for \textsc{Dominating Set}, \textsc{Feedback Vertex Set}, and \textsc{$k$-Coloring} for every fixed~$k$ when a~2-sequence is given.

We further show that, for every fixed~$t \ge 2$, every graph of twin-width at most~$2$ with no induced~$K_{t,t}$ has tree-independence number at most $84(t-1)$.
Given a~$2$-sequence of an~$n$-vertex graph in this class, we construct such a~tree decomposition in time~$n^{O(t)}$.
Finally, we prove that \textsc{Minimum Coloring} is NP-complete on graphs of twin-width at most~$2$, even when a~$2$-sequence is given.
\end{abstract}

\section{Introduction}

A~contraction sequence of a~graph~$G$ starts with the partition of $V(G)$ into singletons and successively merges pairs of parts until a~single part remains.
At each step, two parts are joined by a~red edge whenever there is both an edge and a~nonedge of~$G$ between them.
A~contraction sequence in which every part has red degree at most~$d$ at every step is called a~\emph{$d$-sequence}.
The \emph{twin-width} of~$G$, introduced by Bonnet, Kim, Thomass\'{e}, and Watrigant~\cite{BKTW22-TwwI}, is the least~$d$ for which~$G$ admits a~$d$-sequence.

Graphs of twin-width at most~$1$ have bounded clique-width~\cite{BKRT22-TwwVI}.
Thus many problems, including \textsc{Dominating Set}, \textsc{Feedback Vertex Set}, and \textsc{$k$-Coloring} for every fixed~$k$, can be solved in polynomial time on this class~\cite{CourcelleMR00}.
On the other hand, it was observed that most NP-hard graph problems remain intractable on graphs of twin-width at~most~4 (see for instance~\cite{B26-Coloring}).
Indeed, twin-width is at~most~$4$ for grid graphs (the finite induced subgraphs of the infinite square grid)~\cite[Lemma~9]{BKRTW21-Kernels} and for long subdivisions (obtained by subdividing each edge of an~$n$-vertex graph at least~$2\log n$ times)~\cite[Theorem~2]{BBD22-Tww4}.
Does a~contraction sequence of width~$2$ or~$3$ still give enough structure to solve problems such as \textsc{Dominating Set} in polynomial time?

Clique-width is already unbounded on graphs of twin-width at most~$2$~\cite{BGKTW24-TwwIII,GR00}.
Given a~$2$-sequence, \textsc{Max Independent Set} (and thus \textsc{Max Clique} and \textsc{Vertex Cover}) can still be solved in polynomial time~\cite{BGKTW24-TwwIII}.
On graphs of twin-width at most~$3$, \textsc{Dominating Set}, \textsc{Feedback Vertex Set}, and~\textsc{$3$-Coloring} are NP-complete, even when a~$3$-sequence is given~\cite{BC26-MIS}.
For these last three problems, the case of twin-width~$2$ remained open even when a~$2$-sequence was given.

Vatshelle~\cite{Vatshelle12} introduced the \emph{mim-width} graph parameter.
The mim-width of a~branch decomposition is the maximum, over its cuts, of the size of a~largest induced matching in the bipartite graph formed by the crossing edges.
The mim-width of a~graph is the minimum of this quantity over its branch decompositions.
Many graph problems admit polynomial-time algorithms when a~branch decomposition of bounded mim-width is given~\cite{BJL23}.

\paragraph{Our contributions.}
Our main contribution implies that many problems can be solved in polynomial time on graphs of twin-width at~most~$2$ when a~$2$-sequence is given.
We indeed show that a~$2$-sequence can be turned into a~branch decomposition of bounded mim-width.

\begin{theorem}\label{thm:main}
Every graph~$G$ of twin-width at most~$2$ admits a~branch decomposition of mim-width at most~$8$.
Given a~$2$-sequence of~$G$, one can construct such a~decomposition in polynomial time.
\end{theorem}

At each step of a~$2$-sequence, the red graph is a~disjoint union of paths and cycles.
We introduce the \emph{fringes of red edges} and use them to construct an auxiliary partition at each step.
For each part of this partition, the edges to its complement can be covered by at most~$4$ chain subgraphs.
Each part of the next partition is the union of at most~$5$ parts of the current one.
These partitions therefore define a~rooted tree of bounded degree, which we turn into a~branch decomposition of mim-width at most~$8$.

The decomposition from Theorem~\ref{thm:main} allows us to apply known algorithms for graphs of bounded mim-width~\cite{BJL23,BPT22,JKT20:mim1,JKT20:mim2}.

\begin{corollary}
On graphs given with a~$2$-sequence, \textsc{$k$-Coloring} can be solved in polynomial time for every fixed~$k$.
The same holds for \textsc{Minimum Dominating Set}, \textsc{Total Dominating Set}, \textsc{Independent Dominating Set}, \textsc{Connected Dominating Set}, \textsc{Weighted Feedback Vertex Set}, \textsc{Subset Feedback Vertex Set}, and \textsc{Longest Induced Path}.
\end{corollary}

In the decision version of \textsc{Minimum Coloring}, the number~$k$ of available colors is part of the input.
Graphs of twin-width at most~$1$ are permutation graphs~\cite{AJKPRW25-TwwOne}, on which \textsc{Minimum Coloring} can be solved in polynomial time~\cite{EPL72-Permutation}.
Bonnet~\cite{B26-Coloring} proved that the problem is NP-complete on graphs of twin-width at most~$3$.
In Section~5 we lower this bound to~$2$, even when a~$2$-sequence is given.
Table~\ref{tab:low-twin-width-complexity} summarizes the complexity of these and related problems at low twin-width.

\begin{table}[!htb]
\centering
\small
\begingroup
\setlength{\tabcolsep}{7pt}
\renewcommand{\arraystretch}{1.15}
\begin{tabular}{lcccc}
\toprule
Problem & $\mathcal{T}_1$ & $\mathcal{T}_2$ & $\mathcal{T}_3$ & $\mathcal{T}_4$ \\
\midrule
\textsc{Minimum Coloring}
  & $\mathrm{P}$~\cite{AJKPRW25-TwwOne,EPL72-Permutation}
  & \textbf{NP-c}
  & NP-c~\cite{B26-Coloring}
  & NP-c~\cite{B26-Coloring} \\
\textsc{$3$-Coloring}
  & $\mathrm{P}$~\cite{BKRT22-TwwVI,CourcelleMR00}
  & $\boldsymbol{\mathrm{P}^{\ast}}$
  & NP-c~\cite{BC26-MIS}
  & NP-c~\cite{B26-Coloring} \\
\textsc{$k$-Coloring} (fixed~$k\ge4$)
  & $\mathrm{P}$~\cite{BKRT22-TwwVI,CourcelleMR00}
  & $\boldsymbol{\mathrm{P}^{\ast}}$
  & NP-c~\cite{BC26-MIS}
  & NP-c~\cite{B26-Coloring} \\
\textsc{Max Independent Set}
  & $\mathrm{P}$~\cite{BKRT22-TwwVI,CourcelleMR00}
  & $\mathrm{P}^{\ast}$~\cite{BGKTW24-TwwIII}
  & ?
  & NP-c~\cite{BBD22-Tww4} \\
\textsc{Vertex Cover}
  & $\mathrm{P}$~\cite{BKRT22-TwwVI,CourcelleMR00}
  & $\mathrm{P}^{\ast}$~\cite{BGKTW24-TwwIII}
  & ?
  & NP-c~\cite{BBD22-Tww4} \\
\textsc{Max Clique}
  & $\mathrm{P}$~\cite{BKRT22-TwwVI,CourcelleMR00}
  & $\mathrm{P}^{\ast}$~\cite{BGKTW24-TwwIII}
  & ?
  & NP-c~\cite{BBD22-Tww4} \\
\textsc{Feedback Vertex Set}
  & $\mathrm{P}$~\cite{BKRT22-TwwVI,CourcelleMR00}
  & $\boldsymbol{\mathrm{P}^{\ast}}$
  & NP-c~\cite{BC26-MIS}
  & NP-c~\cite{BBD22-Tww4} \\
\textsc{Dominating Set}
  & $\mathrm{P}$~\cite{BKRT22-TwwVI,CourcelleMR00}
  & $\boldsymbol{\mathrm{P}^{\ast}}$
  & NP-c~\cite{BC26-MIS}
  & NP-c~\cite{BBD22-Tww4} \\
\textsc{Max Induced Matching}
  & $\mathrm{P}$~\cite{BKRT22-TwwVI,CourcelleMR00}
  & $\mathrm{P}^{\ast}$~\cite{BCDGK26-RCML}
  & ?
  & NP-c~\cite{BBD22-Tww4} \\
\textsc{Hamiltonian Path}
  & $\mathrm{P}$~\cite{AJKPRW25-TwwOne,DS94-Hamiltonian}
  & ?
  & ?
  & NP-c~\cite{BBD22-Tww4,IPS82-GridHamiltonian} \\
\textsc{Hamiltonian Cycle}
  & $\mathrm{P}$~\cite{AJKPRW25-TwwOne,DS94-Hamiltonian}
  & ?
  & ?
  & NP-c~\cite{BBD22-Tww4,IPS82-GridHamiltonian} \\
\bottomrule
\end{tabular}
\endgroup
\caption{Complexity of the decision versions on the class~$\mathcal{T}_d$ of graphs of twin-width at most~$d$.
Bold entries are results of this paper.
An asterisk indicates that the polynomial-time algorithm takes a~$d$-sequence as part of the input.}
\label{tab:low-twin-width-complexity}
\end{table}

Table~\ref{tab:twin-width-two-widths} summarizes the width bounds for classes of twin-width at most~$2$.
A~graph class is \emph{weakly sparse} if it excludes~$K_{t,t}$ as a~subgraph for some fixed~$t\ge2$.
A~class is \emph{$K_{t,t}$-free} if it excludes~$K_{t,t}$ as an induced subgraph.
A~class is \emph{edge-stable} if there is an integer~$s\ge1$ such that no graph in the class contains~$2s$ distinct vertices~$a_1,\ldots,a_s,b_1,\ldots,b_s$ with~$a_i$ adjacent to~$b_j$ precisely when~$i\le j$.
Adjacencies among the~$a_i$ and among the~$b_j$ are unrestricted.

Bergougnoux et al.~\cite[Theorem~1.2]{BGGHPS23} showed that every weakly sparse class of twin-width at most~$2$ has bounded treewidth.
Heinrich et al.~\cite[Theorem~26]{HLPR26} showed that every edge-stable class of twin-width at most~$2$ has bounded rank-width and hence bounded clique-width~\cite[Proposition~6.3]{Oum06}.
Theorem~\ref{thm:main} and its proof recover both boundedness results.
The treewidth bound follows directly because bounded mim-width implies bounded treewidth on weakly sparse classes~\cite[Corollary~2]{BMPY25}.
Since complementation preserves the red graphs of a~$2$-sequence of~$G$, the branch decomposition we construct from these red graphs has mim-width at most~$8$ for both~$G$ and its complement.
By Ding et al.~\cite[Corollary~2.4]{DOOV96}, a~cut of sufficiently large rank in an edge-stable class has an induced matching of size~$9$ in its bipartite graph or in its bipartite complement.
Hence every cut of this decomposition has rank bounded in terms of~$s$.

The \emph{tree-independence number} of~$G$, denoted by $\mathrm{tree}\text{-}\alpha(G)$, is the least~$k$ such that~$G$ has a~tree decomposition in which every bag has independence number at most~$k$.
For every fixed~$t\ge2$, Theorem~\ref{thm:main} and a~result of Brettell et al.~\cite[Theorem~6]{BMPY25} imply that~$K_{t,t}$-free graphs of twin-width at most~$2$ have bounded tree-independence number.
In Section~\ref{sec:tree-alpha} we prove the explicit bound~$84(t-1)$.
Given a~$2$-sequence of such an~$n$-vertex graph, we construct a~tree decomposition witnessing this bound in time~$n^{O(t)}$.

\begin{table}[!htb]
\centering
\small
\begingroup
\setlength{\tabcolsep}{8pt}
\renewcommand{\arraystretch}{1.25}
\begin{tabular}{ll}
\toprule
Additional assumption & Bounded parameter \\
\midrule
None & Mim-width \\
Edge-stability & Clique-width \\
Weak sparsity & Treewidth \\
$K_{t,t}$-freeness for fixed~$t\ge2$ & Tree-independence number \\
\bottomrule
\end{tabular}
\endgroup
\caption{Width bounds for graph classes of twin-width at most~$2$.}
\label{tab:twin-width-two-widths}
\end{table}

\paragraph{Organization.}
Section~2 fixes the notation and recalls the relevant width parameters.
In Section~3 we construct the branch decomposition and derive its algorithmic consequences.
In Section~\ref{sec:tree-alpha} we bound the tree-independence number, and in Section~5 we prove the NP-completeness of \textsc{Minimum Coloring}.

\section{Preliminaries}
All graphs in this paper are finite and simple.
For a~positive integer~$n$, let $[n] = \{1, \ldots, n\}$.
For a~graph~$G$, we denote its vertex set and edge set by $V(G)$ and $E(G)$, respectively.
For~$F\subseteq E(G)$, we write~$G-F$ for the graph with vertex set $V(G)$ and edge set $E(G)\setminus F$.
For~$X \subseteq V(G)$, let~$\overline{X} = V(G) \setminus X$.
We denote by $G[X]$ the subgraph induced by~$X$.

For two disjoint sets~$X,Y\subseteq V(G)$, we denote by $G[X,Y]$ the bipartite graph with vertex set~$X\cup Y$ whose edge set consists of all edges of~$G$ with one endpoint in~$X$ and the other in~$Y$. We say that~$X$ is \emph{complete} to~$Y$ if every vertex of~$X$ is adjacent to every vertex of~$Y$, and \emph{anti-complete} to~$Y$ if there is no edge between~$X$ and~$Y$. We say that~$X$ and~$Y$ are \emph{homogeneous} if they are either complete or anti-complete to each other.

Given a~graph~$G$ and a~positive integer~$k$, the decision version of \textsc{Minimum Coloring} asks whether~$G$ admits a~proper coloring with at most~$k$ colors, that is, whether~$\chi(G)\le k$.

\subsection{Induced matchings and mim-width}

A~set~$M$ of edges of a~graph~$G$ is an \emph{induced matching} if the subgraph induced by the endpoints of the edges in~$M$ consists exactly of the edges in~$M$.

For two disjoint sets~$X,Y\subseteq V(G)$, let $\mim_G(X,Y)$ denote the maximum size of an induced matching in $G[X,Y]$.

A~bipartite graph is a~\emph{chain graph} if the neighborhoods of the vertices on one side are linearly ordered by inclusion.
Equivalently, a~bipartite graph is a~chain graph if it has no induced matching of size~$2$.
In particular, if $G[X,Y]$ is a~chain graph, then $\mim_G(X,Y)\le 1$.
A~\emph{chain subgraph} of $G[X,Y]$ is a~chain graph obtained from $G[X,Y]$ by deleting vertices and edges.
A~family of chain subgraphs \emph{covers} $G[X,Y]$ if the union of their edge sets is $E(G[X,Y])$.

A~\emph{branch decomposition} of a~graph~$G$ is a~pair $(T,\eta)$, where~$T$ is a~tree of maximum degree at most~$3$ and~$\eta:V(G)\to L(T)$ is a~bijection from $V(G)$ to the set of leaves of~$T$. For each~$e\in E(T)$, let $(A_e,B_e)$ be the bipartition of $V(G)$ corresponding to the two sets of leaves of the components of~$T-e$.

The \emph{mim-width} of a~branch decomposition $(T,\eta)$ is the maximum size of an induced matching in $G[A_e,B_e]$ over all~$e\in E(T)$. The \emph{mim-width} of~$G$, denoted by $\mim(G)$, is the minimum mim-width over all branch decompositions of~$G$. We set $\mim(G)=0$ when~$G$ is an edgeless graph.

\subsection{Contraction sequence}
A~\emph{trigraph} has two disjoint sets of edges, called \emph{black edges} and \emph{red edges}. For a~trigraph~$H$, let $\mathcal{R}(H)$ be the graph with vertex set $V(H)$ and the red edges of~$H$. A~\emph{red component} is a~connected component of $\mathcal{R}(H)$, and the \emph{red degree} of a~vertex is its degree in $\mathcal{R}(H)$. A~\emph{red path} or a~\emph{red cycle} in~$H$ is a~path or a~cycle in $\mathcal{R}(H)$, respectively.

For a~partition~$\mathcal{P}$ of $V(G)$, the \emph{quotient graph}~$G/\mathcal{P}$ is the trigraph with vertex set~$\mathcal{P}$ in which distinct parts~$P, Q$ are connected by a~black edge if they are complete to each other, by a~red edge if they are not homogeneous, and by no edge if they are anti-complete.

For a~non-negative integer~$d$, a~\emph{$d$-sequence} of an~$n$-vertex graph~$G$ is a~sequence~$\mathcal{P}_n,\ldots,\mathcal{P}_1$ of partitions of $V(G)$ satisfying the following:
\begin{itemize}
    \item $\mathcal{P}_n$ is the partition into singletons and~$\mathcal{P}_1=\{V(G)\}$;
    \item for every~$i\in\{2,\ldots,n\}$, the partition~$\mathcal{P}_{i-1}$ is obtained from~$\mathcal{P}_i$ by merging two parts;
    \item for every~$i\in[n]$, each vertex of~$G/\mathcal{P}_i$
    has red degree at most~$d$.
\end{itemize}
The \emph{twin-width} of~$G$, denoted by $\tww(G)$, is the minimum integer~$d$ such that~$G$ admits a~$d$-sequence.

\section{Graphs of twin-width at most~2 have bounded mim-width}

Throughout this section, fix a~$2$-sequence~$\mathcal{P}_n,\ldots,\mathcal{P}_1$ of an~$n$-vertex graph~$G$.
In the next subsection we define \emph{fringes}, vertex sets associated with certain red edges of the contraction sequence.
We use these fringes to build auxiliary partitions in which every part~$X$ satisfies $\mim_G(X,\overline{X})\le4$.
From these partitions we construct a~branch decomposition.

\subsection{Fringes}
Let~$H_i$ be the spanning subgraph of~$G$ obtained by deleting every edge whose endpoints lie in distinct parts of~$\mathcal{P}_i$ joined by a~black edge in~$G/\mathcal{P}_i$.
Thus, if~$P$ and~$Q$ are distinct parts and $PQ$ is not a~red edge in~$G/\mathcal{P}_i$, then~$P$ and~$Q$ are anti-complete in~$H_i$.
In particular, distinct red components are anti-complete in~$H_i$.

We first recall the following result from Heinrich et al.~\cite{HLPR26}.

\begin{lemma}[{\cite[Lemma~21]{HLPR26}}]\label{lem:1-sequence}
    Let~$\mathcal{P}_n, \ldots, \mathcal{P}_1$ be a~$1$-sequence of a~graph~$G$. For every~$i \in [n]$ and every pair of distinct parts~$P, Q \in \mathcal{P}_i$, the graph $G[P, Q]$ is a~chain graph.
\end{lemma}

\begin{proof}
    We proceed by induction along the sequence. For~$\mathcal{P}_n$, every pair of parts induces a~chain graph. Suppose that~$P_1$ and~$P_2$ are contracted into~$P^*$, and let~$Q$ be another part. Since~$Q$ has red degree at most~$1$ before the contraction, at least one of~$P_1$ and~$P_2$ is homogeneous to~$Q$. By induction, the graph between~$Q$ and the other part is a~chain graph. Adding vertices that are complete or anti-complete to~$Q$ preserves the inclusion order of the neighborhoods of the vertices in~$Q$. Hence, $G[P^*, Q]$ is a~chain graph. All other pairs remain unchanged, which completes the proof.
\end{proof}

For each~$i$, let~$\mathcal{E}_i \subseteq E(\mathcal{R}(G/ \mathcal{P}_i))$ consist of the following edges:
\begin{enumerate}[(1)]
    \item every red edge whose endpoints both have red degree~$2$ and whose red component is a~path;
    \item every edge in a~red component that is a~cycle of length at least~$5$.
\end{enumerate}
Thus, for every $BC \in \mathcal{E}_i$, there are four distinct consecutive parts~$A, B, C, D$ of a~path or cycle, and $AD$ is not a~red edge.

Fix $BC\in \mathcal{E}_i$. Let~$\ell>i$ be the smallest index such that the red edges between parts contained in~$B$ and parts contained in~$C$ do not consist of a~single edge $B'C'$ whose endpoints both have red degree~$2$ at~$\mathcal{P}_\ell$. Since~$\mathcal{P}_n$ has no red edges, such an index exists.

We first show that there is at most~$1$ red edge $B'C'$ with~$B'\subseteq B$ and~$C'\subseteq C$ at~$\mathcal{P}_\ell$. Consider the split from~$\mathcal{P}_{\ell-1}$ to~$\mathcal{P}_\ell$. If the split does not involve either endpoint of the unique red edge, no new red edge between~$B$ and~$C$ appears, since homogeneous pairs remain homogeneous. If an endpoint of the unique red edge is split, the other endpoint still has another red neighbor. Since its red degree is at most~$2$, it cannot have red edges to both parts obtained from the split.

We define $\beta_i(BC)$ as follows. If there is no red edge $B'C'$ in~$G/\mathcal{P}_\ell$ with~$B'\subseteq B$ and~$C'\subseteq C$, let $\beta_i(BC)=\emptyset$. Otherwise, let $B'C'$ be the unique such red edge. Exactly one of~$B'$ and~$C'$ has red degree~$1$, and we define $\beta_i(BC)$ to be this part. We call $\beta_i(BC)$ the \emph{fringe} of $BC$.

\begin{figure}[tbp]
\centering
\begin{tikzpicture}[font=\small,
  part/.style={draw,thick,rounded corners=4pt},
  rededge/.style={red!75!black,very thick},
  oldcut/.style={gray!75,densely dotted,very thick},
  newcut/.style={thick,dashed,rounded corners=5pt}]
  \draw[rededge] (-5.4,0)--(-2.4,0)--(1.9,0)--(5.45,0);
  \draw[part,fill=blue!10] (-6.6,-0.6) rectangle (-4.2,0.6);
  \draw[part,fill=blue!10] (-3.2,-0.6) rectangle (-0.8,0.6);
  \draw[part,fill=red!10] (0.5,-0.6) rectangle (3.4,0.6);
  \draw[part,fill=blue!18] (0.62,-0.47) rectangle (1.5,0.47);
  \draw[part,fill=red!10] (4.3,-0.6) rectangle (6.6,0.6);
  \node at (-5.4,0) {$A$};
  \node at (-2,0) {$B$};
  \node at (1.06,0) {$Y$};
  \node at (2.45,0) {$C\setminus Y$};
  \node at (5.45,0) {$D$};
  \node at (1.95,0.98) {$C$};

  \draw[newcut] (-0.8,1.1)--(-0.8,0.79)--(1.66,0.79)--(1.66,-0.79)--(-0.8,-0.79)--(-0.8,-1.1);

  \node at (-3.75,-1.17) {$L=A\cup B\cup Y$};
  \node at (4.1,-1.17) {$R=(C\setminus Y)\cup D$};
  \node at (0, 1.6) {$\emptyset \neq Y=\beta_i(BC)\subseteq C$};

\end{tikzpicture}
\caption{The cut $(L,R)$ in Lemma~\ref{lem:separateLR} when $\varnothing\ne Y=\beta_i(BC)\subseteq C$. At~$\mathcal{P}_\ell$, there are no red edges between~$Y$ and the parts contained in~$R$.}
\label{fig:separteLR}
\end{figure}
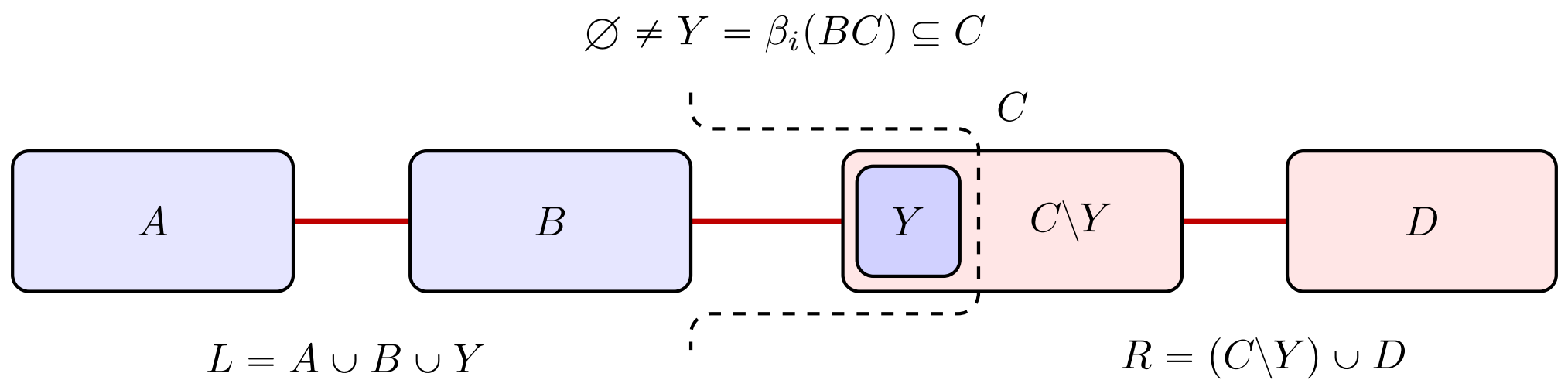

\begin{lemma}\label{lem:separateLR}
Let~$\mathcal{P}_n, \ldots, \mathcal{P}_1$ be a~$2$-sequence of a~graph~$G$. Let~$A,B,C,D$ be four distinct consecutive parts in a~red path or a~red cycle of length at least~$5$ in~$G/\mathcal{P}_i$. Let~$Y=\beta_i(BC)$, and define
\[(L,R)=
\begin{cases}
(A\cup B,\ C\cup D),&Y=\varnothing,\\
(A\cup B\cup Y,\ (C\setminus Y)\cup D),&\emptyset \neq Y\subseteq C,\\
(A\cup(B\setminus Y),\ Y\cup C\cup D),&\emptyset \neq Y\subseteq B.
\end{cases}\]
Then $\mim_{H_i}(L,R)\le 2$.
\end{lemma}

\begin{proof}
    We first consider~$Y=\emptyset$. Consider $G[B,C]$ along the sequence from~$\mathcal{P}_n$ to~$\mathcal{P}_i$. There is no red edge between a~part contained in~$B$ and a~part contained in~$C$ from~$\mathcal{P}_n$ to~$\mathcal{P}_\ell$, and there is at most~$1$ red edge between a~part contained in~$B$ and a~part contained in~$C$ from~$\mathcal{P}_{\ell-1}$ to~$\mathcal{P}_i$. After deleting repeated partitions and then contracting~$B$ and~$C$, we obtain a~$1$-sequence of $G[B,C]$. By Lemma~\ref{lem:1-sequence}, $G[B,C]$ is a~chain graph. Since $AC$, $AD$, and $BD$ are not red in~$G/\mathcal{P}_i$, the only edges of $H_i[L,R]$ are between~$B$ and~$C$. Therefore, $\mim_{H_i}(L,R)\leq 1$.

    We next consider $\emptyset \neq Y\subseteq C$.
    At~$\mathcal{P}_\ell$, the unique red neighbor of~$Y$ is contained in~$B$.
    Hence, there is no red edge between a~part contained in~$B$ and a~part contained in~$C\setminus Y$.
    The graph $G[B,C\setminus Y]$ is a~chain graph, trivially if~$Y=C$ and by the preceding argument otherwise.
    
    Since $AC$, $AD$, and $BD$ are not red edges in~$G/\mathcal{P}_i$,~$A$ is anti-complete to~$C$ and~$D$, and~$B$ is anti-complete to~$D$ in~$H_i$. Also, every vertex of~$Y$ has the same neighborhood in $(C\setminus Y)\cup D$. Indeed, $(C\setminus Y)\cup D$ is a~union of parts of~$\mathcal{P}_\ell$, and none of these parts is a~red neighbor of~$Y$, since the unique red neighbor of~$Y$ is contained in~$B$. Thus, every part contained in $(C\setminus Y)\cup D$ is homogeneous to~$Y$.
    
    The edges between~$Y$ and~$C\setminus Y$ are in~$H_i$ since both are contained in~$C$, and the edges between~$Y$ and~$D$ are in~$H_i$ since $CD$ is a~red edge in~$G/\mathcal{P}_i$. Thus, $H_i[A\cup B \cup Y,(C\setminus Y)\cup D]$ consists of the chain graph between~$B$ and~$C\setminus Y$ and a~complete bipartite graph between~$Y$ and a~subset of $(C\setminus Y)\cup D$. An induced matching contains at most~$1$ edge from each of these two graphs. Therefore, $\mim_{H_i}(A\cup B \cup Y,(C\setminus Y)\cup D)\leq 2$.

    The case $\emptyset\neq Y\subseteq B$ is symmetric.
\end{proof}

\begin{lemma}\label{lem:disjointfringe}
    Let~$\mathcal{P}_n, \ldots, \mathcal{P}_1$ be a~$2$-sequence of a~graph~$G$. Let~$A, B, C$ be distinct parts such that $AB, BC \in \mathcal{E}_i$. If $\beta_i(AB) \subseteq B$ and $\beta_i(BC) \subseteq B$, then $\beta_i(AB) \cap \beta_i(BC) = \emptyset$.
\end{lemma}

\begin{proof}
    Let~$Y = \beta_i(AB)$ and~$Z = \beta_i(BC)$. Suppose that~$Y \cap Z \neq \emptyset$. Since~$Y$ and~$Z$ are parts appearing in the sequence, we obtain either~$Y\subseteq Z$ or~$Z\subseteq Y$.

    Suppose first that~$Y\subsetneq Z$. Let~$\mathcal{P}_j$ and~$\mathcal{P}_k$ be the partitions where~$Y$ and~$Z$ are chosen, respectively. Then~$j>k$. When~$Z$ is chosen at~$\mathcal{P}_k$, its unique red neighbor is contained in~$C$, so~$Z$ is homogeneous to every part contained in~$A$. Since subsets of homogeneous parts are homogeneous,~$Y$ has no red neighbor contained in~$A$ at~$\mathcal{P}_j$, a~contradiction. The case~$Z\subsetneq Y$ is symmetric.

    Suppose now that~$Y=Z$. Let~$\mathcal{P}_j$ and~$\mathcal{P}_k$ be the partitions where~$Y$ is chosen for $AB$ and $BC$, respectively. Suppose first that~$j>k$. When~$Y$ is chosen for $BC$ at~$\mathcal{P}_k$, its unique red neighbor is contained in~$C$, so it is homogeneous to every part contained in~$A$. Since subsets of homogeneous parts are homogeneous,~$Y$ cannot have a~red neighbor contained in~$A$ at~$\mathcal{P}_j$. However, when~$Y$ is chosen for $AB$ at~$\mathcal{P}_j$, its unique red neighbor is contained in~$A$, a~contradiction. The case~$k >  j$ is symmetric. Therefore,~$Y\cap Z=\emptyset$.
\end{proof}

\subsection{Auxiliary partitions}
We define an auxiliary partition~$\mathcal{A}_i$ of $V(G)$ for each~$i\in[n]$. For~$v\in V(G)$, let~$P$ be the part of~$\mathcal{P}_i$ containing~$v$. If~$v\in \beta_i(PQ)\subseteq P$ for some edge $PQ\in\mathcal{E}_i$, set $\lambda_i(v)=Q$. Otherwise, set $\lambda_i(v) = P$. By Lemma~\ref{lem:disjointfringe}, there is at most~$1$ such edge $PQ$, so~$\lambda_i$ is well-defined. Each vertex is assigned only once, according to the part of~$\mathcal{P}_i$ that contains it.

For every connected component~$C$ of $\mathcal{R}(G / \mathcal{P}_i) - \mathcal{E}_i$, let $X_i(C)=\{v\in V(G):\lambda_i(v)\in V(C)\}$.
Define
\[\mathcal{A}_i=\{X_i(C):C\text{ is a~connected component of }\mathcal{R}(G / \mathcal{P}_i)-\mathcal{E}_i \,\text{ and }\, X_i(C)\neq\emptyset\}.\]
The vertex sets of the connected components of $\mathcal{R}(G/\mathcal{P}_i)-\mathcal{E}_i$ form a~partition of~$\mathcal{P}_i$. Since $\lambda_i(v)\in\mathcal{P}_i$ for every~$v\in V(G)$, each vertex~$v$ belongs to exactly one set $X_i(C)$. Therefore, the nonempty sets $X_i(C)$ form a~partition~$\mathcal{A}_i$ of $V(G)$. The parts of~$\mathcal{A}_i$ are not necessarily unions of parts of~$\mathcal{P}_i$.

If a~red component is a~path on at least~$4$ parts, then after deleting~$\mathcal{E}_i$, its two end edges form two components, and every remaining internal part forms a~singleton component. In particular, a~red path on~$4$ parts is divided by deleting its middle edge. A~red path on at most~$3$ parts remains one component. If a~red component is a~cycle of length at least~$5$, every part becomes a~singleton component after deleting~$\mathcal{E}_i$. A~red~$C_3$ or~$C_4$ remains a~single component.

\begin{observation}\label{obs}
Every part of~$\mathcal{A}_i$ meets at most~$3$ parts of~$\mathcal{P}_i$, except for a~part corresponding to a~red~$C_4$, which meets~$4$ parts. Every part of~$\mathcal{P}_i$ meets at most~$3$ parts of~$\mathcal{A}_i$.
\end{observation}

\begin{proof}
Suppose first that a~component of $\mathcal{R}(G/\mathcal{P}_i)-\mathcal{E}_i$ consists of a~single internal part. The corresponding set in~$\mathcal{A}_i$ contains vertices from this part and may also receive vertices from its at most~$2$ red neighbors. Hence, it meets at most~$3$ parts of~$\mathcal{P}_i$.

Suppose next that a~component consists of an endpoint of the red path and its neighbor. Only one edge of~$\mathcal{E}_i$ is incident with one of these two parts and another part outside the component. Thus, the corresponding part of~$\mathcal{A}_i$ meets at most~$1$ part besides these two.

A~red path on at most~$3$ parts, or a~red~$C_3$ or~$C_4$, has no edge in~$\mathcal{E}_i$. Thus, the corresponding part of~$\mathcal{A}_i$ meets only the parts in that red component. Hence, every part of~$\mathcal{A}_i$ meets at most~$3$ parts of~$\mathcal{P}_i$, except for a~part corresponding to a~red~$C_4$, which meets~$4$ parts.

We now show that every part of~$\mathcal{P}_i$ meets at most~$3$ parts of~$\mathcal{A}_i$. Let~$P\in\mathcal{P}_i$. The vertices of~$P$ that do not belong to any $\beta_i(PQ)$ belong to one part of~$\mathcal{A}_i$. Since at most~$2$ edges of~$\mathcal{E}_i$ are incident with~$P$, the vertices in the sets $\beta_i(PQ)$ belong to at most~$2$ other parts of~$\mathcal{A}_i$. Therefore,~$P$ meets at most~$3$ parts of~$\mathcal{A}_i$.
\end{proof}

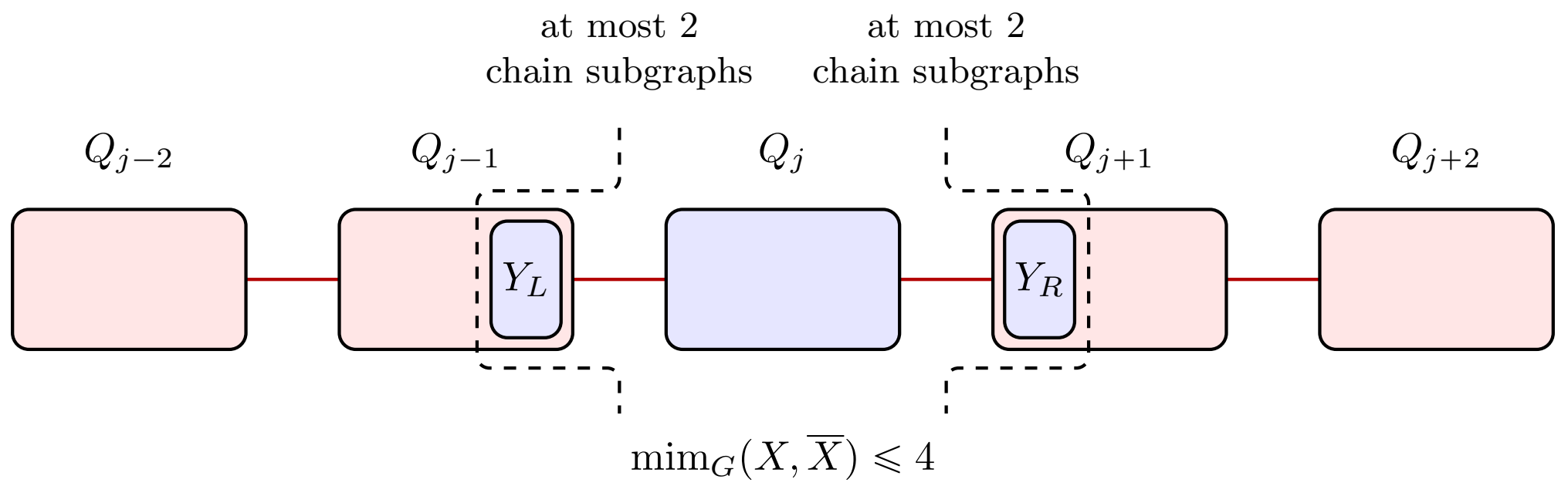
\begin{figure}[tbp]
    \centering
    \begin{tikzpicture}[font=\small,
      part/.style={draw,thick,rounded corners=4pt},
      cut/.style={thick,dashed,rounded corners=3pt}]
      \draw[red!70!black,thick] (-5.6,0)--(5.6,0);
      \foreach \x in {-5.6,-2.8,2.8,5.6} {
        \draw[part,fill=red!10] ({\x-1},-.6) rectangle ({\x+1},.6);
      }
      \draw[part,fill=blue!10] (-1,-.6) rectangle (1,.6);
      \draw[part,fill=blue!10] (-2.5,-.5) rectangle (-1.9,.5);
      \draw[part,fill=blue!10] (1.9,-.5) rectangle (2.5,.5);
      \node[inner sep=0pt] at (-2.2,0) {$Y_L$};
      \node[inner sep=0pt] at (2.2,0) {$Y_R$};
      \node at (-5.6,1.08) {$Q_{j-2}$};
      \node at (-2.8,1.08) {$Q_{j-1}$};
      \node at (0,1.08) {$Q_j$};
      \node at (2.8,1.08) {$Q_{j+1}$};
      \node at (5.6,1.08) {$Q_{j+2}$};
      \draw[cut] (-1.4,1.3)--(-1.4,.76)--(-2.62,.76) --(-2.62,-.76)--(-1.4,-.76)--(-1.4,-1.15);
      \draw[cut] (1.4, 1.3)--(1.4,.76)--(2.62,.76) --(2.62,-.76)--(1.4,-.76)--(1.4,-1.15);
      \node[above,align=center,font=\footnotesize,inner sep=2pt] at (-1.4,1.55)
        {at most~$2$\\chain subgraphs};
      \node[above,align=center,font=\footnotesize,inner sep=2pt] at (1.4,1.55)
        {at most~$2$\\chain subgraphs};
      \node at (0,-1.52) {$\mim_G(X,\overline X)\le4$};
    \end{tikzpicture}
\caption{The auxiliary part~$X=Y_L\cup Q_j\cup Y_R$ is shown in blue, where~$Y_L$ and~$Y_R$ are the fringes assigned to~$Q_j$. For each dashed cut, at most~$2$ chain subgraphs cover the crossing edges within parts or between parts connected by a~red edge. The remaining crossing edges are between parts that are complete to each other. Adding these edges preserves the chain property.}
\label{fig:mim-seven}
\end{figure}

\begin{lemma}\label{lem:Aimimwidthbound}
For every~$i\in[n]$ and~$X\in\mathcal{A}_i$, we have $\mim_G(X,\overline X)\le4$.
\end{lemma}
\begin{proof}
The proof of Lemma~\ref{lem:separateLR} shows that the edges of $H_i[L,R]$ for each local cut $(L,R)$ are covered by a~chain graph and at most~$1$ complete bipartite graph involving the fringe.
Every complete bipartite graph is a~chain graph.

Suppose first that~$X$ corresponds to a~singleton component consisting of an internal part~$Q_j$ of a~red path or a~part~$Q_j$ of a~red cycle of length at least~$5$. Let~$Q_1,\ldots,Q_t$ be the parts of the red component in order if it is a~red path, and in cyclic order if it is a~red cycle. By the construction of~$\mathcal{A}_i$,~$X\subseteq Q_{j-1}\cup Q_j\cup Q_{j+1}$.

Consider the two cuts from Lemma~\ref{lem:separateLR} for $Q_{j-1}Q_j$ and $Q_jQ_{j+1}$. For each of the two cuts,~$X$ lies entirely on one side. Every neighbor of~$X$ in~$H_i$ is contained in one of~$Q_{j-2},\ldots,Q_{j+2}$. If such a~vertex is not in~$X$, then it lies on the side not containing~$X$ in at least one of the two cuts. For a~red cycle, the indices are taken modulo its length.

Let~$X_j=X\cap Q_j$.
By construction,~$X$ consists of~$X_j$ and at most~$2$ fringes assigned to~$Q_j$ from its red neighbors.
The chain graph from each cut has its side in~$X$ contained in~$X_j$.
For each fringe assigned to~$Q_j$, the complete bipartite graph from its local cut has that entire fringe as its side in~$X$.
Thus, at most~$4$ chain subgraphs cover the crossing edges in~$H_i$.

We now add the crossing edges not in~$H_i$ to these subgraphs.
For each part~$Q\in\mathcal P_i$ meeting~$X$, these edges incident with~$X\cap Q$ form a~complete bipartite graph, since they join~$X\cap Q$ to vertices in parts that are complete to~$Q$.
Add all such edges incident with~$X_j$ to one of the two chain graphs, first adding any missing vertices of~$X_j$ as isolated vertices.
These edges add the same neighbor set to every vertex of~$X_j$, so their neighborhoods remain nested.
Every other part~$Q$ meeting~$X$ contributes a~fringe~$Y=X\cap Q$ assigned to~$Q_j$.
Add its crossing edges not in~$H_i$ to the complete bipartite graph whose side in~$X$ is~$Y$.
This graph remains complete bipartite.
Hence, the same four subgraphs cover all edges of $G[X,\overline X]$.

Suppose next that~$X$ corresponds to a~component consisting of an endpoint of a~red path on at least~$4$ parts and its neighbor.
There is only one edge of~$\mathcal E_i$ between this component and the rest of the red path, so one local cut suffices.
We add the crossing edges not in~$H_i$ incident with the neighbor and any incoming fringe to the two subgraphs in the same way.
One additional complete bipartite graph covers those incident with the endpoint.
Thus, at most~$3$ chain subgraphs suffice.
If~$X$ corresponds to a~red path on at most~$3$ parts or to a~red~$C_3$ or~$C_4$, there are no crossing edges in~$H_i$, and at most~$4$ complete bipartite graphs cover the remaining edges.
An induced matching contains at most~$1$ edge from each chain subgraph. Therefore, $\mim_G(X,\overline X)\le 4$. This completes the proof.
\end{proof}

\subsection{Contraction and lift}

Fix~$i\in\{2,\ldots,n\}$. Suppose that two parts~$P_1,P_2\in\mathcal{P}_i$ are contracted into~$P^*\in\mathcal{P}_{i-1}$. Let~$f:\nobreak\mathcal{P}_i\to\nobreak\mathcal{P}_{i-1}$ map each part to the part containing it. Then, we have $f(P_1)=f(P_2)=P^*$, and $f(P)=P$ for every~$P\in\mathcal{P}_i\setminus \{P_1,P_2\}$. We say that an edge $PQ\in E(\mathcal{R}(G / \mathcal{P}_i))$ is a~\emph{red lift} of an edge $BC\in E(\mathcal{R}(G / \mathcal{P}_{i-1}))$ if $\{f(P),f(Q)\}=\{B,C\}$.

\begin{lemma}\label{lem:lift}
Let~$e\in\mathcal{E}_{i-1}$. Then, exactly one of the following holds.
\begin{enumerate}[{\rm (1)}]
\item The edge~$e$ has no red lift and $\beta_{i-1}(e)=\emptyset$.
\item The edge~$e$ has exactly one red lift $PQ\in\mathcal{E}_i$ and $\beta_{i-1}(e)=\beta_i(PQ)$. If $\emptyset \neq \beta_i(PQ)\subseteq P$, then $\lambda_{i-1}(v) = f(Q)$ for every~$v \in \beta_i(PQ)$, and if $\emptyset \neq \beta_i(PQ)\subseteq Q$, then $\lambda_{i-1}(v) = f(P)$ for every~$v \in \beta_i(PQ)$.
\item The edge~$e$ has exactly one red lift $PQ$, where~$P$ has red degree~$1$ and~$Q$ has red degree~$2$ in~$G/\mathcal{P}_i$. Then, $\beta_{i-1}(e)=P$ and $\lambda_{i-1}(v) = f(Q)$ for every~$v \in P$.
\end{enumerate}
\end{lemma}

\begin{proof}
We first show that~$e$ has at most~$1$ red lift. If neither endpoint of~$e$ is~$P^*$, then the only possible lift is the same edge. Suppose that one endpoint is~$P^*$, and let~$R$ be the other endpoint. Since~$e\in\mathcal{E}_{i-1}$, the part~$R$ has another red neighbor different from~$P^*$. This part is unchanged at~$\mathcal{P}_i$. Thus, at most~$1$ of $RP_1$ and $RP_2$ is a~red edge at~$\mathcal{P}_i$, since otherwise~$R$ would have red degree at least~$3$. Hence,~$e$ has at most~$1$ red lift.

Suppose first that~$e$ has no red lift. Then, there is no red edge between a~part contained in one endpoint of~$e$ and a~part contained in the other endpoint at~$\mathcal{P}_i$. Hence, $\beta_{i-1}(e)=\emptyset$, and (1) holds.

Suppose next that the red lift is $PQ$ and both~$P$ and~$Q$ have red degree~$2$ at~$\mathcal{P}_i$. Suppose that $f(P) =  P^*$, and let~$P' \in \mathcal{P}_i$ be the part such that $f(P') = P^*$ with~$P' \neq P$. Since $PQ$ is the unique red lift of~$e$,~$P'$ is homogeneous to~$Q$. Since subsets of homogeneous parts are homogeneous, these parts remain homogeneous at every~$\mathcal{P}_j$ with~$j>i$. Thus, the unique red edge between the two parts is the same for~$e$ and $PQ$ at every~$\mathcal{P}_j$ with~$j \ge i$. Hence, $\beta_{i-1}(e)=\beta_i(PQ)$.

We show that $PQ\in\mathcal{E}_i$. Since both~$P$ and~$Q$ have red degree~$2$, the only case in which $PQ\notin\mathcal{E}_i$ is when $PQ$ belongs to a~red~$C_3$ or~$C_4$. If no two parts of this cycle are contracted together, then they form the same red cycle at~$\mathcal{P}_{i-1}$. So, we have~$e\notin\mathcal{E}_{i-1}$, a~contradiction. If two parts of a~red~$C_3$ are contracted, the remaining part has red degree~$1$. If two adjacent parts of a~red~$C_4$ are contracted, they form a~red~$C_3$ at~$\mathcal{P}_{i-1}$. If two opposite parts are contracted, they form a~red path on~$3$ parts whose endpoints have red degree~$1$. In each case,~$e\notin\mathcal{E}_{i-1}$, a~contradiction. Therefore, $PQ\in\mathcal{E}_i$.

If $\emptyset \neq \beta_i(PQ)\subseteq P$, then $\lambda_i(v) = Q$ for every~$v \in \beta_i(PQ)$. Since $\beta_{i-1}(e)=\beta_i(PQ)$ and the endpoints of~$e$ are $f(P)$ and $f(Q)$, we have $\beta_{i-1}(e) \subseteq f(P)$. Thus, we obtain $\lambda_{i-1}(v) = f(Q)$ for every~$v \in \beta_i(PQ)$. The case $\emptyset \neq \beta_i(PQ)\subseteq Q$ is symmetric. Thus, (2) holds.

It remains to consider the case where an endpoint of the unique red lift has red degree~$1$. Suppose that this endpoint is~$P$. We show that~$Q$ has red degree~$2$. Suppose that~$Q$ also has red degree~$1$. If one of~$P$ and~$Q$ is contracted into~$P^*$, then the other part still has red degree~$1$ at~$\mathcal{P}_{i-1}$, contradicting~$e\in\mathcal{E}_{i-1}$. Suppose that neither~$P$ nor~$Q$ is contracted. For both~$P$ and~$Q$ to have red degree~$2$ at~$\mathcal{P}_{i-1}$, both $PP^*$ and $QP^*$ must be red. Then~$P$,~$Q$, and~$P^*$ form a~red~$C_3$, so~$e\notin\mathcal{E}_{i-1}$, a~contradiction. Hence,~$Q$ has red degree~$2$.

At~$\mathcal{P}_i$, the unique red lift is $PQ$, where~$P$ has red degree~$1$ and~$Q$ has red degree~$2$. Thus, by the definition of $\beta_{i-1}(e)$, we have $\beta_{i-1}(e)=P$. By the definition of~$\lambda_{i-1}$, we have $\lambda_{i-1}(v) = f(Q)$ for every~$v \in \beta_{i-1}(e)$. Hence, (3) holds.
\end{proof}

\begin{lemma}\label{lem:refines}
For every~$X\in\mathcal{A}_i$, there exists~$X'\in\mathcal{A}_{i-1}$ such that~$X\subseteq X'$.
\end{lemma}
\begin{proof}
We compare~$\lambda_i$ and~$\lambda_{i-1}$. For each part~$Q\in\mathcal{P}_i$, let
\[B_i(Q)=\{v\in V(G):\lambda_i(v)=Q\}.\]
The set $B_i(Q)$ can be empty.

\begin{claim}\label{claim:receiver}
Suppose that~$Q$ has red degree~$2$ in~$G/\mathcal{P}_i$ and its red component is neither a~$C_3$ nor a~$C_4$. Then, for every~$v \in B_i(Q)$, the parts $\lambda_{i-1}(v)$ and $f(Q)$ belong to the same connected component of $\mathcal{R}(G / \mathcal{P}_{i-1})-\mathcal{E}_{i-1}$.
\end{claim}

\begin{clproof}
Let~$v\in B_i(Q)$. Suppose first that~$v\in Q$. Since $\lambda_i(v)=Q$, the vertex~$v$ does not belong to any $\beta_i(QR)\subseteq Q$ with $QR\in\mathcal{E}_i$. Suppose that~$v \in \beta_{i-1}(e)$ for some~$e \in \mathcal{E}_{i-1}$. By Lemma~\ref{lem:lift}, the nonempty set $\beta_{i-1}(e)$ is either equal to $\beta_i(PR)$ for some edge $PR\in\mathcal{E}_i$ or a~part of~$\mathcal{P}_i$ of red degree~$1$. In the first case,~$v$ does not belong to $\beta_i(PR)$, and in the second case, the fringe cannot contain~$v$ since~$Q$ has red degree~$2$. So, we have~$v \notin \beta_{i-1}(e)$ for every~$e \in \mathcal{E}_{i-1}$. Hence, $\lambda_{i-1}(v)=f(Q)$.

We now suppose that~$v\notin Q$. Then,~$v\in\beta_i(PQ)\subseteq P$ for some $PQ\in\mathcal{E}_i$. In particular, both~$P$ and~$Q$ have red degree~$2$ in~$G/\mathcal{P}_i$. If $f(P)f(Q)\in\mathcal{E}_{i-1}$, then by Lemma~\ref{lem:lift}, we have $\beta_{i-1}(f(P)f(Q))=\beta_i(PQ)$, and hence $\lambda_{i-1}(v)=f(Q)$. Suppose that $f(P)f(Q)\notin\mathcal{E}_{i-1}$. If $f(P)=f(Q)$, then $f(P)$ and $f(Q)$ belong to the same connected component of $\mathcal{R}(G/\mathcal{P}_{i-1})-\mathcal{E}_{i-1}$. Suppose now that $f(P) \neq f(Q)$. Then, we have $f(P)f(Q) \in E(\mathcal{R}(G/\mathcal{P}_{i-1}))$. As $f(P)f(Q)\notin\mathcal{E}_{i-1}$, $f(P)$ and $f(Q)$ belong to the same component of $\mathcal{R}(G/\mathcal{P}_{i-1})-\mathcal{E}_{i-1}$.

It remains to determine $\lambda_{i-1}(v)$. By Lemma~\ref{lem:lift}, every nonempty fringe at~$\mathcal{P}_{i-1}$ is either equal to $\beta_i(PR)$ for some $PR\in\mathcal{E}_i$ or a~part of~$\mathcal{P}_i$ of red degree~$1$. By Lemma~\ref{lem:disjointfringe},~$v$ does not belong to any other fringe in the first case. Moreover,~$v$ cannot belong to a~fringe in the second case, since~$v\in P$ and~$P$ has red degree~$2$. Hence,~$v \notin \beta_{i-1}(e)$ for every~$e \in \mathcal{E}_{i-1}$. Thus, $\lambda_{i-1}(v)=f(P)$. Since $f(P)$ and $f(Q)$ are in the same connected component, the claim follows.
\end{clproof}

Let~$C$ be a~connected component of $\mathcal{R}(G/\mathcal{P}_i)-\mathcal{E}_i$. We show that all vertices of $X_i(C)$ belong to one part of~$\mathcal{A}_{i-1}$. Suppose first that~$C$ consists of one part~$Q$ of red degree~$2$. Then $X_i(C)=B_i(Q)$, and the result follows from Claim~\ref{claim:receiver}.

Now, we suppose that~$C$ consists of an endpoint~$P$ of a~red path and its neighbor~$Q$, where~$P$ has red degree~$1$ and~$Q$ has red degree~$2$. Since no edge of~$\mathcal{E}_i$ is incident with~$P$, we obtain $X_i(C)=\nobreak P\cup B_i(Q)$. By Claim~\ref{claim:receiver}, all vertices of $B_i(Q)$ are assigned by~$\lambda_{i-1}$ to parts in the connected component containing $f(Q)$. If $f(P)=f(Q)$, then they are the same part. Suppose that $f(P)\neq f(Q)$. If $f(P)f(Q)\notin\mathcal{E}_{i-1}$, then $f(P)f(Q)$ is a~red edge, so $f(P)$ and $f(Q)$ belong to the same connected component of $\mathcal{R}(G/\mathcal{P}_{i-1})-\mathcal{E}_{i-1}$. If $f(P)f(Q)\in\mathcal{E}_{i-1}$, then by Lemma~\ref{lem:lift} (3), we have $\beta_{i-1}(f(P)f(Q))=P$, and hence $\lambda_{i-1}(v) = f(Q)$ for every~$v \in P$. Since~$Q$ is the only red neighbor of~$P$ at~$\mathcal{P}_i$, Lemma~\ref{lem:lift} implies that if $\beta_{i-1}(e)$ contains a~vertex of~$P$ for some~$e\in\nobreak\mathcal{E}_{i-1}$, then~$e=f(P)f(Q)$. Hence, no other fringe $\beta_{i-1}(e)$ contains a~vertex of~$P$. Thus, there is a~connected component~$C'$ of $\mathcal R(G/\mathcal{P}_{i-1})-\mathcal{E}_{i-1}$ such that $\lambda_{i-1}(v)\in V(C')$ for every~$v\in X_i(C)$. Hence, $X_i(C)\subseteq X_{i-1}(C')$, and therefore $X_i(C)$ is contained in one part of~$\mathcal A_{i-1}$. The same argument applies when~$C$ is a~red path on~$3$ parts.

Now we suppose that~$C$ consists of an isolated part~$P$. By Lemma~\ref{lem:lift}, no nonempty set $\beta_{i-1}(e)$ contains a~vertex of~$P$. Hence, $\lambda_{i-1}(v)=f(P)$ for every~$v\in P$.

Suppose that~$C$ is a~red component on~$2$ parts~$P$ and~$Q$. By the proof of Lemma~\ref{lem:lift}, the edge between $f(P)$ and $f(Q)$ does not belong to~$\mathcal{E}_{i-1}$ whenever $f(P)\neq f(Q)$. Also, no other red lift is incident with~$P$ or~$Q$, since both~$P$ and~$Q$ have red degree~$1$. Hence, either $f(P)=f(Q)$, or $f(P)f(Q)$ is a~red edge of~$G/\mathcal{P}_{i-1}$ that does not belong to~$\mathcal{E}_{i-1}$. Therefore, all vertices of $X_i(C)$ belong to one part of~$\mathcal{A}_{i-1}$.

Finally, suppose that~$C$ is a~red~$C_3$ or~$C_4$. By Lemma~\ref{lem:lift}, for every edge $PQ$ of~$C$, either $f(P)=f(Q)$ or $f(P)f(Q)\notin\mathcal{E}_{i-1}$. Also, no set $\beta_{i-1}(e)$ contains a~vertex of~$C$. Indeed, no edge of~$C$ belongs to~$\mathcal{E}_i$, and every part of~$C$ has red degree~$2$. Thus, for every part~$P$ of~$C$ and every~$v\in P$, we have $\lambda_{i-1}(v)=f(P)$. By Lemma~\ref{lem:lift}, for every edge $PQ$ of~$C$, either $f(P)=f(Q)$ or $f(P)f(Q)\notin\mathcal{E}_{i-1}$. Hence, all parts $f(P)$ with~$P\in V(C)$ belong to the same connected component of $\mathcal{R}(G/\mathcal{P}_{i-1})-\mathcal{E}_{i-1}$.

These cases cover every connected component of $\mathcal{R}(G/\mathcal{P}_i)-\mathcal{E}_i$. Therefore, every part of~$\mathcal{A}_i$ is contained in a~part of~$\mathcal{A}_{i-1}$. This completes the proof.
\end{proof}

\begin{lemma}\label{lem:15bound}
Every part of~$\mathcal{A}_{i-1}$ is the union of at most~$5$ parts of~$\mathcal{A}_i$.
\end{lemma}

\begin{proof}
Fix a~contraction from~$\mathcal P_i$ to~$\mathcal P_{i-1}$, and let~$f:\mathcal P_i\to\mathcal P_{i-1}$ be the map used in Lemma~\ref{lem:refines}.
For each component~$C$ of $\mathcal R(G/\mathcal P_i)-\mathcal E_i$, we choose a~part~$Q_C\in V(C)$ such that every vertex of $X_i(C)$ is assigned by~$\lambda_{i-1}$ to the component containing $f(Q_C)$.

If~$C$ consists of one part~$Q$ of red degree~$2$, choose~$Q_C=Q$; Claim~\ref{claim:receiver} gives the required property.
If~$C$ consists of an endpoint~$P$ and its neighbor~$Q$ of red degree~$2$, again choose~$Q_C=Q$.
Claim~\ref{claim:receiver} applies to every vertex assigned to~$Q$ by~$\lambda_i$, and the endpoint argument in the proof of Lemma~\ref{lem:refines} shows that every vertex of~$P$ is assigned to the component containing $f(Q)$.
For a~red path on~$3$ parts, choose its middle part and apply the same argument to its two endpoints.
In all remaining cases,~$C$ is an isolated part, a~red component on~$2$ parts, or a~red~$C_3$ or~$C_4$.
The corresponding cases in the proof of Lemma~\ref{lem:refines} show that all parts $f(P)$ for~$P\in V(C)$ belong to one component and all vertices of $X_i(C)$ are assigned to that component.
Thus, any part of~$C$ can be chosen as~$Q_C$.

Let~$D$ be a~component of $\mathcal R(G/\mathcal P_{i-1})-\mathcal E_{i-1}$ with $X_{i-1}(D)\ne\emptyset$.
If a~nonempty old part $X_i(C)$ is contained in $X_{i-1}(D)$, then $f(Q_C)\in V(D)$, by the property above.
Distinct components~$C$ have distinct chosen parts~$Q_C$.
Consequently, the number of old parts contained in $X_{i-1}(D)$ is at most $|f^{-1}(V(D))|$.
Every component~$D$ has at most~$4$ parts, and only the contracted part can have~$2$ preimages under~$f$.
Hence, $|f^{-1}(V(D))|\le |V(D)|+1\le 5$.
Every part of~$\mathcal A_{i-1}$ is therefore the union of at most~$5$ parts of~$\mathcal A_i$, which proves the result.
\end{proof}

\subsection{Constructing the branch decomposition}\label{sec:decomposition}

We first construct the rooted tree used in the proof of Theorem~\ref{thm:main}.

\begin{lemma}\label{lem:tree}
Given a~graph~$G$ and a~$2$-sequence of~$G$, we can construct in polynomial time a~rooted tree~$T$ whose leaves are the vertices of~$G$ such that every internal node has at most~$5$ children and $\mim_G(X_t,\overline{X_t})\le 4$ for every node~$t$ of~$T$, where~$X_t$ is the set of leaves in the subtree rooted at~$t$.
\end{lemma}

\begin{proof}
By construction,~$\mathcal{A}_n$ consists of singletons and~$\mathcal{A}_1=\{V(G)\}$. We start with one leaf for each part of~$\mathcal{A}_n$. For each~$i=n,n-1,\ldots,2$ and each~$X\in\mathcal{A}_{i-1}$, if~$X\in\mathcal{A}_i$, we keep the node for~$X$. Otherwise, we create a~node for~$X$ and make the nodes corresponding to the parts of~$\mathcal{A}_i$ contained in~$X$ its children.

By Lemma~\ref{lem:refines}, every part of~$\mathcal{A}_{i-1}$ is a~union of parts of~$\mathcal{A}_i$. If~$X\in\mathcal{A}_{i-1}$ does not belong to~$\mathcal{A}_i$, then~$X$ contains at least~$2$ parts of~$\mathcal{A}_i$. Hence, every internal node has at least~$2$ children. Since~$\mathcal{A}_n$ consists of~$n$ singletons, the tree has~$n$ leaves, and therefore at most~$n-1$ internal nodes. Thus, the tree has at most $2n-1$ nodes.

By Lemma~\ref{lem:15bound}, every part of~$\mathcal{A}_{i-1}$ is the union of at most~$5$ parts of~$\mathcal{A}_i$. Hence, every internal node has at most~$5$ children. For every node~$t$, the set~$X_t$ belongs to some~$\mathcal{A}_i$. Thus, by Lemma~\ref{lem:Aimimwidthbound}, $\mim_G(X_t,\overline{X_t})\leq 4$.

All sets used in the construction can be computed in polynomial time. We first compute $\mathcal{R}(G/\mathcal{P}_i)$ and~$\mathcal{E}_i$ for every~$i\in[n]$. For each edge in~$\mathcal{E}_i$, we consider at most~$n$ partitions to determine its fringe. Since each~$\mathcal{E}_i$ has at most~$n$ edges and there are at most~$n$ partitions, all fringes can be computed in polynomial time. We can then compute the maps~$\lambda_i$ and the auxiliary partitions~$\mathcal{A}_i$. Finally, comparing~$\mathcal{A}_i$ and~$\mathcal{A}_{i-1}$ for every~$i$ gives the tree. Thus, the tree can be constructed in polynomial time.
\end{proof}

\begin{proof}[Proof of Theorem~\ref{thm:main}]
Let~$G$ be a~graph of twin-width at most~$2$, and let~$\mathcal{P}_n,\ldots,\mathcal{P}_1$ be a~$2$-sequence of~$G$.
The empty graph and the graph on~$1$ vertex have mim-width~$0$, so we may assume that $|V(G)|\ge 2$.
Let~$T$ be the rooted tree given by Lemma~\ref{lem:tree}.
Every internal node has at most~$5$ children, so the underlying unrooted tree has maximum degree at most~$6$.

Let~$t$ be an internal node of~$T$.
The sets of leaves in the components of~$T-t$ partition $V(G)$.
Each such set~$X$ satisfies $\mim_G(X,\overline X)\le 4$ by Lemma~\ref{lem:tree}.
Indeed, it is either the set of leaves below a~child of~$t$, or the complement of the set of leaves below~$t$.
The two sides of a~cut have the same induced matching bound.

For each internal node~$t$ of degree at least~$4$, partition its incident edges into three groups, each of size at most~$2$.
Replace~$t$ by a~central node.
For each group of two edges, introduce one node adjacent to the central node and attach the two edges of the group to this new node.
Attach every edge in a~group of size~$1$ directly to the central node.
We leave the nodes of degree at most~$3$ unchanged.
Let~$T'$ be the resulting tree.
Every node of~$T'$ has degree at most~$3$, and its leaves are the same as those of~$T$.

Every edge of~$T'$ corresponding to an original edge of~$T$ defines the same cut, and hence has mim-width at most~$4$.
Every other edge joins a~central node to the node for a~group of two edges.
One side of its cut is the union of the sets of leaves in two components of~$T-t$ for the corresponding original node~$t$.

Let~$X_1,\ldots,X_s$ be the sets of leaves in at most~$2$ of these components, where~$s\le 2$, and let~$X=\bigcup_{j=1}^s X_j$.
Let~$M$ be an induced matching in $G[X,\overline X]$.
For each~$j\in[s]$, let~$M_j$ be the set of edges of~$M$ with an endpoint in~$X_j$.
Then~$M_j$ is an induced matching in $G[X_j,\overline{X_j}]$.
By the bound above, $|M_j|\le 4$.
Therefore,
\[
|M|=\sum_{j=1}^s |M_j|\le 4s\le 8.
\]
Thus, every edge cut of~$T'$ has mim-width at most~$8$.

Let~$\eta$ be the bijection between $V(G)$ and the leaves of~$T'$ inherited from~$T$.
Since~$T'$ has maximum degree at most~$3$, $(T',\eta)$ is a~branch decomposition of~$G$ with mim-width at most~$8$.

By Lemma~\ref{lem:tree}, the rooted tree~$T$ can be constructed in polynomial time from the given $2$-sequence.
Each internal node is replaced by a~tree with at most~$4$ nodes, so~$T'$ and the branch decomposition can also be constructed in polynomial time.
This completes the proof.
\end{proof}

\subsection{Graph problems}
\begin{corollary}
Let~$G$ be a~graph of twin-width at most~$2$, and suppose that a~witnessing contraction sequence is given. Then, the following problems can be solved in polynomial time on~$G$\textup{:} \textsc{$k$-Coloring} for every fixed~$k$, \textsc{Minimum Dominating Set}, \textsc{Total Dominating Set}, \textsc{Independent Dominating Set}, \textsc{Connected Dominating Set}, \textsc{Weighted Feedback Vertex Set}, \textsc{Subset Feedback Vertex Set} and \textsc{Longest Induced Path}.
\end{corollary}

\begin{proof}
    Bergougnoux, Dreier, and Jaffke~\cite{BJL23} showed that, given a~branch decomposition of mim-width at most~$w$, every problem expressible in their $\mathsf{A \&C~DN}$ logic can be solved in time~$n^{O(w)}$. This includes Minimum Dominating Set, Total Dominating Set, Independent Dominating Set, Connected Dominating Set, and~$k$-Coloring for every fixed~$k$.

    Jaffke, Kwon, and Telle showed that Longest Induced Path can be solved in time~$n^{O(w)}$~\cite{JKT20:mim1}, and that Weighted Feedback Vertex Set can be solved in time~$n^{O(w)}$~\cite{JKT20:mim2}, when a~branch decomposition of mim-width at most~$w$ is given. Bergougnoux, Papadopoulos, and Telle~\cite{BPT22} showed that Subset Feedback Vertex Set can be solved in time~$n^{O(w^2)}$ under the same assumption.

    By Theorem~\ref{thm:main}, from the given $2$-sequence of~$G$, we can construct in polynomial time a~branch decomposition of~$G$ with mim-width at most~$8$. Thus, by taking~$w=8$ in the algorithms above, all the listed problems can be solved in polynomial time on~$G$.
\end{proof}

\section{\texorpdfstring{$K_{t,t}$-free}{K(t,t)-free} graphs of twin-width at most~2 have bounded \texorpdfstring{tree-$\alpha$}{tree-alpha}}\label{sec:tree-alpha}

In this section, we consider~$K_{t,t}$-free graphs for a~fixed integer~$t\ge 2$, where~$K_{t,t}$ is excluded as an induced subgraph.
For~$S\subseteq V(G)$, let $\alpha(S)$ denote the maximum size of an independent set in $G[S]$, and let $N_G(S)$ be the set of vertices with a~neighbor in~$S$.
A~\emph{tree decomposition} of~$G$ is a~tree~$T$ together with bags~$B_x\subseteq V(G)$ for~$x\in V(T)$ such that every vertex belongs to a~bag, every edge has both endpoints in a~bag, and the nodes whose bags contain a~fixed vertex induce a~connected subtree of~$T$.
The independence number of this decomposition is $\max_{x\in V(T)}\alpha(B_x)$.
The \emph{tree-independence number} of~$G$, denoted by $\mathrm{tree}\text{-}\alpha(G)$, is the minimum of these numbers over all tree decompositions of~$G$.

\begin{theorem}\label{thm:tree-independence}
Let~$t\ge 2$ and let~$G$ be a~$K_{t,t}$-free graph of twin-width at most~$2$. Then
\[
\mathrm{tree}\text{-}\alpha(G)\le 84(t-1).
\]
Given a~$2$-sequence of~$G$, we can construct a~tree decomposition of~$G$ with independence number at most $84(t-1)$ in time $|V(G)|^{O(t)}$.
\end{theorem}

We use the rooted tree constructed in Section~\ref{sec:decomposition}, before replacing its internal nodes by trees of maximum degree~$3$.

\begin{lemma}\label{lem:chain-cuts-tree-independence}
Let~$t\ge 2$ and let~$r,d\ge 1$ be integers.
Let~$G$ be a~nonempty~$K_{t,t}$-free graph, and let~$T$ be a~rooted tree whose leaves are the vertices of~$G$ and whose internal nodes have at most~$d$ children.
For a~node~$x$, let~$X_x$ be the set of leaves in its subtree.
Suppose that every cut $G[X_x,\overline{X_x}]$ is covered by at most~$r$ chain subgraphs with sides in~$X_x$ and $\overline{X_x}$.
Then there is a~tree decomposition of~$G$ whose underlying tree is~$T$ and whose independence number is at most $r(1+dr)(t-1)$.
\end{lemma}

\begin{proof}
Let~$K=r(t-1)$.
For every set~$X=X_x$ corresponding to a~node of~$T$, define
\[
H_X=\{v\in\overline X:\alpha(N_G(v)\cap X)>K\},
\qquad
I_X=N_G(\overline X\setminus H_X)\cap X.
\]
Here,~$H_X$ contains the outside vertices whose neighborhoods in~$X$ have large independence number, and~$I_X$ contains the endpoints in~$X$ of the remaining cut edges.
The set~$C_X=H_X\cup I_X$ covers every edge of $G[X,\overline X]$.
Indeed, if $uv$ is such an edge with~$u\in X$ and~$v\in\overline X$, then either~$v\in H_X$ or~$u\in I_X$.

We first show that
\begin{equation}\label{eq:chain-cut-independence-bounds}
\alpha(H_X)\le K
\qquad\text{and}\qquad
\alpha(I_X)\le rK.
\end{equation}
Let~$F_1,\ldots,F_r$ be chain subgraphs whose edge sets cover the cut; add empty graphs if necessary.
We add isolated vertices so that each~$F_j$ has bipartition $(X,\overline X)$.
For each~$v\in H_X$, at least one neighborhood $N_{F_j}(v)\cap X$ contains an independent set of size~$t$.
Otherwise, an independent set in $N_G(v)\cap X$ would have size at most $r(t-1)=K$.
Assign~$v$ to one such index~$j$.
If~$H_X$ contained an independent set of size~$K+1$, at least~$t$ of its vertices would be assigned to the same chain graph.
Their neighborhoods in~$X$ are nested, and the smallest of these neighborhoods contains an independent set of size~$t$.
These two independent sets induce a~$K_{t,t}$, a~contradiction.
Thus, $\alpha(H_X)\le K$.

Let~$L=\overline X\setminus H_X$.
For each~$j$, the neighborhoods in~$X$ of vertices of~$L$ in~$F_j$ are nested.
Their union is therefore the neighborhood of one vertex of~$L$, unless~$L$ is empty.
The union has independence number at most~$K$, by the definition of~$L$.
These at most~$r$ sets cover~$I_X$, so $\alpha(I_X)\le rK$.

We now show three relations between these sets.
If~$Y\subseteq X$ are sets corresponding to nodes of~$T$, then
\begin{equation}\label{eq:chain-cut-nesting}
I_X\cap Y\subseteq I_Y
\qquad\text{and}\qquad
H_Y\setminus X\subseteq H_X.
\end{equation}
For the first relation, a~vertex of~$\overline X\setminus H_X$ also belongs to~$\overline Y\setminus H_Y$.
The second relation follows directly from $N_G(v)\cap Y\subseteq N_G(v)\cap X$.
If~$Y$ and~$Z$ correspond to distinct children of the same node, then
\begin{equation}\label{eq:chain-cut-siblings}
H_Y\cap Z\subseteq I_Z.
\end{equation}
Indeed, let~$v\in H_Y\cap Z$.
Since $\alpha(N_G(v)\cap Y)>K$ and $\alpha(H_Z)\le K$, some neighbor of~$v$ in~$Y$ does not belong to~$H_Z$.
This neighbor lies in~$\overline Z\setminus H_Z$, so~$v\in I_Z$.

For each nonroot node~$x$, label the edge joining~$x$ to its parent by~$C_{X_x}$.
For a~vertex~$v\in V(G)$, let~$T_v$ consist of its leaf and all tree edges whose labels contain~$v$.
We show that~$T_v$ is connected.
Suppose first that~$v\in I_{X_x}$.
By the first relation in~\eqref{eq:chain-cut-nesting}, every edge on the path from~$x$ to the leaf~$v$ has a~label containing~$v$.
Suppose next that~$v\in H_{X_x}$.
Let~$z$ be the lowest common ancestor of~$x$ and the leaf~$v$, and let~$y$ and~$w$ be its children whose subtrees contain~$x$ and the leaf~$v$, respectively.
By the second relation in~\eqref{eq:chain-cut-nesting}, all edges from~$x$ to~$z$ have labels containing~$v$, and~$v\in H_{X_y}$.
By~\eqref{eq:chain-cut-siblings},~$v\in I_{X_w}$.
By the first relation in~\eqref{eq:chain-cut-nesting}, all edges from~$z$ to the leaf~$v$ have labels containing this vertex.
Thus, every edge of~$T_v$ is connected to the corresponding leaf within~$T_v$.

Set~$B_x=\{v:x\in V(T_v)\}$ for every node~$x$ of~$T$.
Every vertex belongs to the bag at its leaf, and the nodes whose bags contain it form the connected subtree~$T_v$.
Let $uv\in E(G)$.
If the leaf~$v$ belongs to~$T_u$, the bag at that leaf contains both endpoints.
Otherwise, consider the first edge $ab$ on the path from the leaf~$u$ to the leaf~$v$ with~$a\in V(T_u)$ and~$b\notin V(T_u)$.
Its label covers a~cut separating~$u$ and~$v$, so it contains~$u$ or~$v$.
It cannot contain~$u$, since then $ab$ would belong to~$T_u$.
Hence, it contains~$v$, and~$a\in V(T_u)\cap V(T_v)$.
The bag~$B_a$ therefore contains both endpoints of $uv$.
Thus, the bags form a~tree decomposition of~$G$.

For an internal node~$x$ with children~$y_1,\ldots,y_s$, where~$s\le d$, relations~\eqref{eq:chain-cut-nesting} and~\eqref{eq:chain-cut-siblings} give
\[
B_x=H_{X_x}\cup\bigcup_{j=1}^s I_{X_{y_j}}.
\]
To see this, note that~$B_x$ is the union of the labels of its incident edges.
The set~$I_{X_x}$ is contained in the union of the sets~$I_{X_{y_j}}$.
Every vertex of~$H_{X_{y_j}}$ either lies outside~$X_x$, and hence belongs to~$H_{X_x}$, or lies in another child set~$X_{y_h}$, and hence belongs to~$I_{X_{y_h}}$.
At the root, take~$C_{V(G)}=\emptyset$.
Consequently, by~\eqref{eq:chain-cut-independence-bounds},
\[
\alpha(B_x)\le K+srK\le r(1+dr)(t-1).
\]
For a~leaf~$x$ corresponding to~$v$, we have~$H_{\{v\}}=\emptyset$, since~$K\ge 1$, and hence~$B_x=\{v\}$.
The same bound therefore holds for every bag.
\end{proof}

\begin{proof}[Proof of Theorem~\ref{thm:tree-independence}]
The result is immediate for the empty graph, so suppose that~$G$ is nonempty.
Take the rooted tree constructed in the proof of Lemma~\ref{lem:tree}.
Its internal nodes have at most~$5$ children, and every set of leaves below a~node is a~part of one of the auxiliary partitions.
By the proof of Lemma~\ref{lem:Aimimwidthbound}, every corresponding cut is covered by at most~$4$ chain subgraphs.
Lemma~\ref{lem:chain-cuts-tree-independence}, with~$r=4$ and~$d=5$, gives
\[
\mathrm{tree}\text{-}\alpha(G)
\le 4(1+5\cdot 4)(t-1)=84(t-1).
\]

Given a~$2$-sequence, we can construct this rooted tree in polynomial time.
For each set~$X$ of leaves below a~node and each~$v\in\overline X$, we check all subsets of $N_G(v)\cap X$ of size $4(t-1)+1$.
The vertex~$v$ belongs to~$H_X$ if and only if at least one of these subsets is independent.
After computing all sets~$H_X$, we can compute all sets~$I_X$ and construct the bags in polynomial time.
Since the tree has $O(|V(G)|)$ nodes, the total running time is $|V(G)|^{O(t)}$. This completes the proof.
\end{proof}

\suppressfloats[t]

\section{Minimum Coloring}

Let~$A=\{u_1,\ldots,u_k\}$ be an ordered clique. We define three operations, each of which introduces a~new ordered copy~$B=\{v_1,\ldots,v_k\}$ of~$A$, disjoint from all previously introduced cliques, and adds edges between~$A$ and~$B$ as follows.

\begin{itemize}
    \item We obtain $\ce(A)$ by adding the edges
    $E(A,B)=\{u_pv_q:1\le p<q\le k\}$.

    \item For each~$i\in[k-1]$, we obtain $\cs_i(A)$ by adding the same edges as in $\ce(A)$, except that $u_iv_{i+1}$ is replaced by $u_iv_i$:
    $E(A,B)=\big(\{u_pv_q:1\le p<q\le k\}\setminus\{u_iv_{i+1}\}\big)\cup\{u_iv_i\}$.

    \item For each~$i\in[k-1]$, we obtain $\cc_i(A)$ by adding the same edges as in $\ce(A)$, except that $u_iv_{i+1}$ is omitted:
    $E(A,B)=\{u_pv_q:1\le p<q\le k\}\setminus\{u_iv_{i+1}\}$.
\end{itemize}
In each case, the notation denotes the newly introduced clique~$B$.

\begin{figure}[htbp]
\centering
\begingroup
\colorlet{coloring1}{blue!70!black}
\colorlet{coloring2}{orange!85!black}
\colorlet{coloring3}{teal!75!black}
\colorlet{coloring4}{violet!80!black}
\colorlet{coloring5}{magenta!85!red}
\newcommand{\DrawColoringGadget}[2]{%
  \begin{scope}[xshift=#1cm]
    \draw[block,fill=blue!10] (-.25,-2.00) rectangle (.25,2.00);
    \draw[block,fill=red!10] (1.45,-2.00) rectangle (1.95,2.00);
    \foreach \p in {1,...,5} {
      \coordinate (u\p) at (0,{.85*(3-\p)});
      \coordinate (v\p) at (1.7,{.85*(3-\p)});
    }
    \foreach \p/\q in {1/2,1/3,1/4,1/5,2/4,2/5,3/4,3/5,4/5}
      \draw[edge] (u\p)--(v\q);
    \ifnum#2=0
      \draw[edge] (u2)--(v3);
      \node at (.85,2.90) {$\ce$};
    \else
      \draw[removededge] (u2)--(v3);
      \node[text=orange!90!black,fill=white,inner sep=.5pt,font=\small]
        at (1.19,.255) {$\times$};
      \ifnum#2=1
        \draw[addededge] (u2)--(v2);
        \node at (.85,2.90) {$\cs_2$};
      \else
        \node at (.85,2.90) {$\cc_2$};
        \node at (2.80,2.35) {or};
      \fi
    \fi
    \foreach \p in {1,...,5} {
      \def\outputcolor{\p}
      \ifnum#2=1
        \ifnum\p=2 \def\outputcolor{3}\fi
        \ifnum\p=3 \def\outputcolor{2}\fi
      \fi
      \node[vertex,fill=coloring\p] at (u\p) {};
      \node[vertex,fill=coloring\outputcolor] at (v\p) {};
      \node[anchor=east,inner sep=1pt] at (-.34,{.85*(3-\p)})
        {$c_{\p}$};
      \node[anchor=west,inner sep=1pt] at (2.04,{.85*(3-\p)})
        {$c_{\outputcolor}$};
      \ifnum#2=2
        \def\alternativecolor{\p}
        \ifnum\p=2 \def\alternativecolor{3}\fi
        \ifnum\p=3 \def\alternativecolor{2}\fi
        \node[vertex,fill=coloring\alternativecolor]
          at (2.80,{.85*(3-\p)}) {};
        \node[anchor=west,inner sep=1pt] at (3.14,{.85*(3-\p)})
          {$c_{\alternativecolor}$};
      \fi
    }
    \node at (0,2.35) {$A$};
    \node at (1.7,2.35) {$B$};
  \end{scope}
}
\begin{tikzpicture}[font=\small,
  edge/.style={black!60,line width=.6pt},
  removededge/.style={orange!90!black,line width=1.4pt,
    dash pattern=on 5pt off 3pt,
    preaction={draw=white,solid,line width=3pt}},
  addededge/.style={blue!85!black,line width=1.5pt,
    preaction={draw=white,solid,line width=3pt}},
  vertex/.style={circle,draw,fill=black,inner sep=2.3pt,outer sep=0pt},
  block/.style={dashed,rounded corners=4pt,thick}]
  \DrawColoringGadget{0}{0}
  \DrawColoringGadget{5}{1}
  \DrawColoringGadget{10}{2}
\end{tikzpicture}
\endgroup
\caption{The operations for~$k=5$.
Gray lines show unchanged edges.
Blue solid lines show added edges, and crossed orange dashed lines show omitted edges.
For $\cc_2$, we show both colorings of~$B$.
We omit edges within each clique.}
\label{fig:coloring-gadgets}
\end{figure}
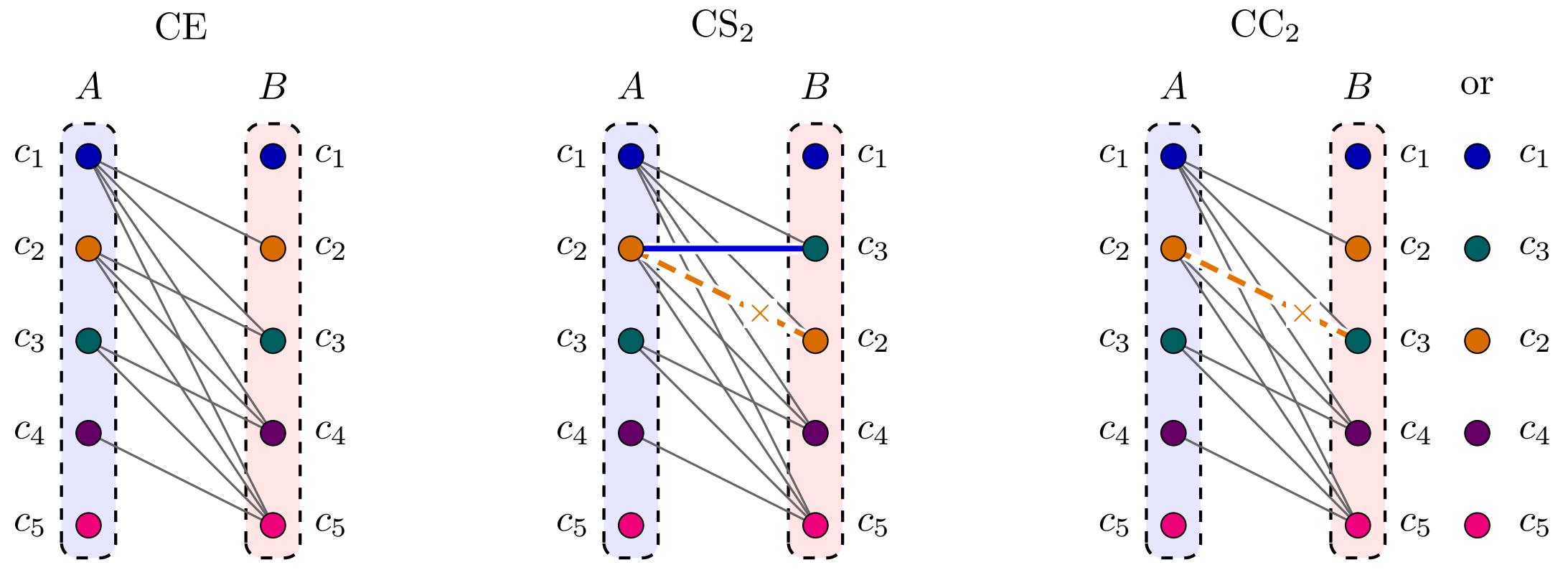

Finally, for $1\le i<j\le k$, we define
\[\cc^*_{(i,j)}(A)=\cs_i\circ\cdots\circ\cs_{j-2}\circ\cc_{j-1}\circ\cs_{j-2}\circ\cdots\circ\cs_i(A),\]
where each operation is applied to the clique introduced by the preceding operation. 

\begin{lemma}\label{lemma:coloringedgegadget}
    Let~$A=\{u_1,\ldots,u_k\}$ be an ordered clique, and let~$f\in\{\ce\}\cup\{\cs_i,\cc_i: i\in[k-1]\}$. Let~$G$ be the graph obtained by applying~$f$ to~$A$, and write~$B=f(A)=\{v_1,\ldots,v_k\}$. Let~$\sigma: V(G) \to [k]$ be a~proper~$k$-coloring of~$G$. Then the following hold.
    \begin{itemize}
        \item If~$f = \ce$, then $\sigma(u_j) = \sigma(v_j)$ for all~$j \in [k]$.
        \item If~$f = \cs_i$, then $\sigma(u_i) = \sigma(v_{i+1})$, $\sigma(u_{i+1}) = \sigma(v_i)$, and $\sigma(u_j) = \sigma(v_j)$ for all~$j \in [k]\setminus \{i, i+1\}$.
        \item If~$f = \cc_i$, then either $\sigma(u_j) = \sigma(v_j)$ for each~$j \in [k]$ or $\sigma(u_i) = \sigma(v_{i+1})$, $\sigma(u_{i+1}) = \sigma(v_i)$, and $\sigma(u_j) = \sigma(v_j)$ for each~$j \in [k]\setminus \{i, i+1\}$.
    \end{itemize}
\end{lemma}
\begin{proof}
    Without loss of generality, we may assume that $\sigma(u_j) = j$ for each~$j \in [k]$.
    
    Suppose first that~$f = \ce$. For each~$j\in[k]$, the vertex~$v_j$ is adjacent to~$u_1,\ldots,u_{j-1}$ and~$v_{j+1},\ldots,v_k$. Considering~$v_k,\ldots,v_1$ in this order, we obtain $\sigma(v_j)=j$ for each~$j\in[k]$.

    Now, we suppose that~$f\in\{\cs_i,\cc_i : i \in [k-1]\}$. The same argument gives $\sigma(v_j)=j$ for each~$j\ge i+2$. Since both~$v_i$ and~$v_{i+1}$ are adjacent to~$u_1, \ldots, u_{i-1}$ and~$v_{i+2}, \ldots, v_k$, we obtain $\{\sigma(v_i),\sigma(v_{i+1})\}=\{i,i+1\}$. Then, considering~$v_{i-1},\ldots,v_1$ in this order, we have $\sigma(v_j)=j$ for each~$j<i$.

    If~$f = \cs_i$, then we obtain $\sigma(v_i) = i+1$ and $\sigma(v_{i+1}) = i$ because $u_iv_i \in E(G)$. If~$f=\cc_i$, then either $\sigma(v_i)=i$ and $\sigma(v_{i+1})=i+1$, or $\sigma(v_i)=i+1$ and $\sigma(v_{i+1})=i$.
\end{proof}

\subsection{NP-completeness}

\begin{theorem}
    \textsc{Minimum Coloring} is NP-complete on graphs of twin-width at most~$2$, even when a~witnessing contraction sequence is given.
\end{theorem}

\begin{proof}
    We reduce from \textsc{Disjoint Paths on Directed Acyclic Graphs}, which is known to be NP-complete~\cite[Theorem~4.8(b)]{GHKLOR14-DigraphWidth}. Let~$\mathcal{D} = (V, E)$ be a~directed acyclic graph and let $(s_1, t_1), \ldots, (s_\ell, t_\ell)$ be pairs of vertices of~$\mathcal{D}$, where~$\ell\ge1$ and all $2\ell$ terminals are distinct. The problem asks whether there exist pairwise vertex-disjoint directed paths~$P_1, \ldots, P_\ell$ in~$\mathcal{D}$ such that~$P_i$ is a~path from~$s_i$ to~$t_i$ for every~$i \in [\ell]$.

    We may assume that each~$s_i$ has no incoming edges and each~$t_i$ has no outgoing edges. Indeed, for each~$i\in[\ell]$, we introduce two new vertices~$s_i'$ and~$t_i'$, add the edges $s_i's_i$ and $t_it_i'$, and replace $(s_i,t_i)$ with $(s_i',t_i')$. This adds $2\ell$ vertices and $2\ell$ edges and preserves both acyclicity and the existence of the required disjoint paths. We use the same notation~$\mathcal{D}=(V,E)$ and $(s_1,t_1),\ldots,(s_\ell,t_\ell)$ for the new instance.

    We now construct a~graph~$G_{\mathcal{D}}$ as follows.
    Let~$V = \{v_1,\ldots,v_n\}$ and~$E = \{e_1,\ldots,e_m\}$, where $v_iv_j\in E$ implies~$i<j$, and let~$k=n+m+1$. We choose a~bijection~$\lambda: V\cup E\to[k-1]$ such that
    \begin{itemize}
        \item $\lambda(s_i)=2i-1$ for each~$i \in [\ell]$,
        \item $\lambda(t_i)=2i$ for each~$i \in [\ell]$, and
        \item $\lambda(e_r)=n+r$ for each~$r\in[m]$.
    \end{itemize}
    We start with an ordered clique~$A_0=\{x_1,\ldots,x_k\}$, where~$x_{\lambda(z)}$ corresponds to~$z$ for each~$z\in V\cup E$. We call the remaining vertex~$x_k$ the \emph{dummy vertex}.
    
    Starting with~$A_0$, we apply the following operations in order.
    
    \begin{enumerate}[(1)]
        \item\label{item:DAGconstruction} We first apply $\ce$ to~$A_0$ and use the newly introduced clique as the starting clique for processing~$v_1$. We consider~$v_1,\ldots,v_n$ in this order. For each~$i\in[n]$, let~$e^-_{i,1},\ldots,e^-_{i,a_i}$ and~$e^+_{i,1},\ldots,e^+_{i,b_i}$ be the incoming and outgoing edges of~$v_i$, respectively, each listed in an arbitrary order.
    
        For~$v_i$, we first successively apply $\cc^*_{(\lambda(v_i),\lambda(e^-_{i,1}))},\ldots,\cc^*_{(\lambda(v_i),\lambda(e^-_{i,a_i}))}$, and let~$A_i$ denote the resulting ordered clique.
        We then successively apply $\cc^*_{(\lambda(v_i),\lambda(e^+_{i,1}))},\ldots,\cc^*_{(\lambda(v_i),\lambda(e^+_{i,b_i}))}$.
    
        The resulting clique is used as the starting clique for the next vertex. If either list of edges is empty, the corresponding sequence of operations is omitted. Thus, each edge~$e=v_iv_j$ is considered twice: first as an outgoing edge of~$v_i$, and then as an incoming edge of~$v_j$.

        Since~$v_1$ has no incoming edges,~$A_1=\ce(A_0)$ is the new clique introduced by the initial operation and is disjoint from~$A_0$.
    
        \item Let~$R=[k-1]\setminus\{2,4,\ldots,2\ell\}$. We successively apply $\cc^*_{(p,q)}$ for all~$p,q\in R$ with~$p<q$, in lexicographic order.
    
        \item We successively apply $\cs_{2h-1}$ for~$h = 1, \ldots, \ell$, in this order.
    \end{enumerate}
    
    Each operation is applied to the clique introduced by the preceding operation, with the first operation applied to~$A_0$. Let~$A_{\text{last}}=\{y_1,\ldots,y_k\}$ be the last ordered clique. Finally, we add the edges
    $E(A_{\text{last}},A_0)=\{y_px_q:1\le p<q\le k\}$.
    
    Let~$G_\mathcal{D}$ be the graph consisting of all cliques and edges introduced in this construction, including all intermediate cliques introduced by $\cc^*$. Thus, the ordered cliques of size~$k$ are arranged in a~cyclic order, with edges between consecutive cliques determined by the operations above.
    
    Let~$H$ be the graph obtained after~(\ref{item:DAGconstruction}). Let~$A_j = \{x^j_1, \ldots, x^j_k\}$ for each~$j \in \{0\} \cup [n]$, where~$x^0_r =\nobreak x_r$ for each~$r \in [k]$.

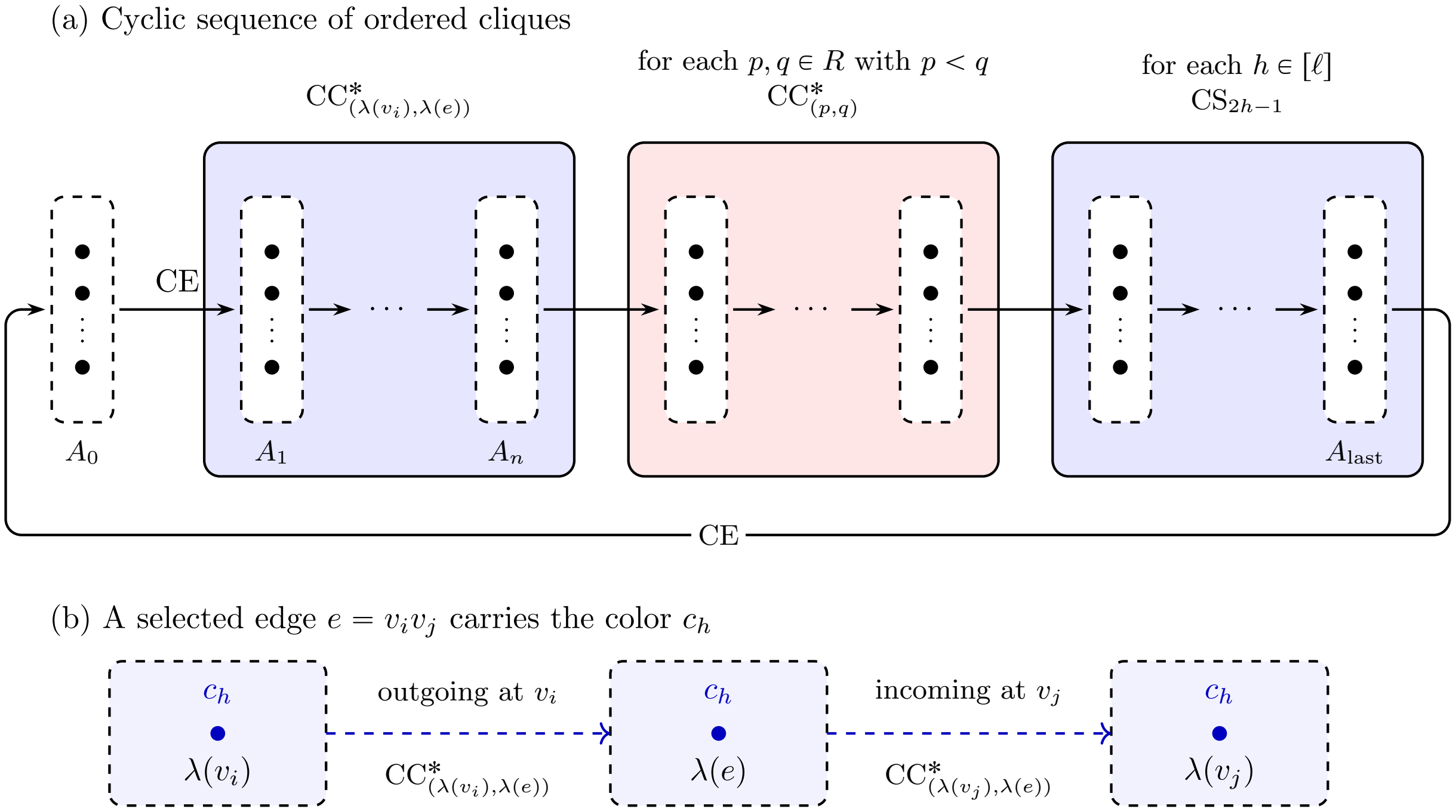
\begin{figure}[tb]
\centering
\begingroup
\newcommand{\ColoringOverviewColumn}[2]{%
  \draw[column,fill=white] (#1-.34,-1.25) rectangle (#1+.34,1.25);
  \foreach \y in {.64,.18,-.64}
    \node[vertex] at (#1,\y) {};
  \foreach \offset in {-.12,0,.12}
    \fill (#1,{(.18-.64)/2+\offset}) circle[radius=.45pt];
  \node[anchor=north,font=\footnotesize,inner sep=0pt] at (#1,-1.45) {#2};
}
\begin{tikzpicture}[
  font=\small,
  vertex/.style={circle,draw,fill=black,inner sep=1.5pt,outer sep=0pt},
  column/.style={thick,dashed,rounded corners=4pt},
  stage/.style={thick,rounded corners=5pt},
  operation/.style={thick,-{Stealth[length=2mm,width=1.35mm]},
    shorten <=2pt,shorten >=2pt},
  flow/.style={blue!75!black,thick,dashed,->},
  position/.style={thick,dashed,rounded corners=4pt,fill=blue!5},
  note/.style={font=\footnotesize,align=center,inner sep=2pt}
]
  \node[anchor=west] at (-.5,3.20) {(a) Cyclic sequence of ordered cliques};

  \draw[stage,fill=blue!10] (1.35,-1.85) rectangle (5.45,1.85);
  \draw[stage,fill=red!10] (6.05,-1.85) rectangle (10.15,1.85);
  \draw[stage,fill=blue!10] (10.75,-1.85) rectangle (14.85,1.85);
  \node[note] at (3.4,2.32)
    {$\cc^*_{(\lambda(v_i),\lambda(e))}$};
  \node[note] at (8.1,2.50)
    {for each~$p,q\in R$ with~$p<q$\\$\cc^*_{(p,q)}$};
  \node[note] at (12.8,2.50)
    {for each~$h\in[\ell]$\\$\cs_{2h-1}$};

  \ColoringOverviewColumn{0}{$A_0$}
  \ColoringOverviewColumn{2.1}{$A_1$}
  \ColoringOverviewColumn{4.7}{$A_n$}
  \ColoringOverviewColumn{6.8}{}
  \ColoringOverviewColumn{9.4}{}
  \ColoringOverviewColumn{11.5}{}
  \ColoringOverviewColumn{14.1}{$A_{\mathrm{last}}$}

  \draw[operation] (.34,0) -- (1.76,0);
  \node at (1.05,.32) {$\ce$};
  \draw[operation] (2.44,0) -- (3.05,0);
  \node at (3.4,0) {$\cdots$};
  \draw[operation] (3.75,0) -- (4.36,0);
  \draw[operation] (5.04,0) -- (6.46,0);
  \draw[operation] (7.14,0) -- (7.75,0);
  \node at (8.1,0) {$\cdots$};
  \draw[operation] (8.45,0) -- (9.06,0);
  \draw[operation] (9.74,0) -- (11.16,0);
  \draw[operation] (11.84,0) -- (12.45,0);
  \node at (12.8,0) {$\cdots$};
  \draw[operation] (13.15,0) -- (13.76,0);

  \draw[operation,rounded corners=5pt]
    (14.44,0) -- (15.15,0) -- (15.15,-2.50)
    -- (-.85,-2.50) -- (-.85,0) -- (-.34,0);
  \node[note,fill=white] at (7.05,-2.50) {$\ce$};

  \node[anchor=west] at (-.5,-3.45)
    {(b) A~selected edge~$e=v_iv_j$ carries the color~$c_h$};
  \foreach \x in {1.5,7.05,12.6} {
    \draw[position] (\x-1.2,-5.50) rectangle (\x+1.2,-3.90);
    \node[vertex,draw=blue!75!black,fill=blue!75!black]
      at (\x,-4.70) {};
    \node[blue!75!black,anchor=south] at (\x,-4.50) {$c_h$};
  }
  \node at (1.5,-5.13) {$\lambda(v_i)$};
  \node at (7.05,-5.13) {$\lambda(e)$};
  \node at (12.6,-5.13) {$\lambda(v_j)$};
  \draw[flow] (2.7,-4.70) -- (5.85,-4.70);
  \draw[flow] (8.25,-4.70) -- (11.4,-4.70);
  \node[note,anchor=south] at (4.275,-4.47) {outgoing at~$v_i$};
  \node[note,anchor=south] at (9.825,-4.47) {incoming at~$v_j$};
  \node[note,anchor=north] at (4.275,-4.93)
    {$\cc^*_{(\lambda(v_i),\lambda(e))}$};
  \node[note,anchor=north] at (9.825,-4.93)
    {$\cc^*_{(\lambda(v_j),\lambda(e))}$};
\end{tikzpicture}
\endgroup
\caption{The construction of~$G_{\mathcal D}$.
In (a), arrows show the operations between ordered cliques.
In (b), the color~$c_h$ moves along~$e=v_iv_j$.
We omit some vertices and all edges within each clique.}
\label{fig:coloring-reduction}
\end{figure}

    \begin{claim}\label{claim:DAGpathpropagate}
        There exist pairwise vertex-disjoint directed paths~$P_1, \ldots, P_\ell$ in~$\mathcal{D}$ such that~$P_i$ is a~path from~$s_i$ to~$t_i$ for each~$i \in [\ell]$ if and only if~$H$ has a~proper~$k$-coloring~$\sigma$ such that $\sigma(x_{\lambda(s_i)}) = \sigma(x^n_{\lambda(t_i)})$ for each~$i \in [\ell]$. Moreover, for each such coloring~$\sigma$, the paths can be chosen so that $V(P_i) = \{v_j \in V(\mathcal{D}) : \sigma(x^j_{\lambda(v_j)}) = c_i\}$ for each~$i \in [\ell]$, where $\sigma(x_{\lambda(s_i)}) = c_i$.
    \end{claim}
    
    \begin{clproof}
        By Lemma~\ref{lemma:coloringedgegadget} and the definition of $\cc^*_{(p,q)}$, this operation either keeps the colors in all positions unchanged or exchanges the colors in positions~$p$ and~$q$. In both cases, we obtain a~proper~$k$-coloring of the resulting graph.
    
        Suppose first that such paths~$P_1,\ldots,P_\ell$ exist. Let~$\sigma$ be a~proper~$k$-coloring of~$A_0$. We first extend~$\sigma$ to~$A_1$ using the initial $\ce$ operation, preserving the color in each position. We consider~$v_1,\ldots,v_n$ in this order. For each vertex~$v_j$, we first consider its incoming edges $wv_j \in E(\mathcal{D})$. If $wv_j$ belongs to one of the paths~$P_1,\ldots,P_\ell$, we extend~$\sigma$ by exchanging the colors in positions $\lambda(v_j)$ and $\lambda(wv_j)$; otherwise, we keep all colors in their positions. We then consider its outgoing edges $v_jw$ in the same way, using positions $\lambda(v_j)$ and $\lambda(v_jw)$. We obtain a~proper~$k$-coloring of~$H$.
    
        Since the paths are vertex-disjoint and each edge $vw$ is considered first at~$v$ and then at~$w$, each color~$c_i$ moves along~$P_i$ and reaches position $\lambda(t_i)$. As~$t_i$ has no outgoing edges, this color remains in position $\lambda(t_i)$ until the end of (\ref{item:DAGconstruction}). Thus, $\sigma(x_{\lambda(s_i)})=\sigma(x^n_{\lambda(t_i)})$ for each~$i\in[\ell]$.
    
        Conversely, suppose that~$H$ has a~$k$-coloring~$\sigma$ such that $\sigma(x_{\lambda(s_i)}) = \sigma(x^n_{\lambda(t_i)})$ for each~$i \in [\ell]$. By the initial $\ce$ operation, $\sigma(x^1_r)=\sigma(x_r)$ for every~$r\in[k]$. For each~$i\in[\ell]$, we follow the color~$c_i$, recording its position in~$A_0$ and in the last clique introduced by each $\cc^*$.
    
        Consider a~vertex~$v_j$ and an incoming edge $v_pv_j\in E(\mathcal{D})$ with~$p<j$. We obtain~$A_j$ after performing $\cc^*_{(\lambda(v_p),\lambda(v_pv_j))}$ and $\cc^*_{(\lambda(v_j),\lambda(v_pv_j))}$ in this order. Hence, $\sigma(x^j_{\lambda(v_pv_j)})=\sigma(x^n_{\lambda(v_pv_j)})$. Since~$A_n$ is a~clique and $\sigma(x^n_{\lambda(t_i)})=c_i$, we obtain $\sigma(x^j_{\lambda(v_pv_j)})\ne c_i$ for each~$i\in[\ell]$. Thus,~$c_i$ cannot move from position $\lambda(v_j)$ to position $\lambda(v_pv_j)$ when this edge is considered at~$v_j$.
    
        Hence,~$c_i$ can leave position $\lambda(v_j)$ only through an outgoing edge $v_jv_q$ with~$j<q$. Then, it must move from position $\lambda(v_jv_q)$ to position $\lambda(v_q)$ when $v_jv_q$ is considered at~$v_q$. Otherwise, we would have $\sigma(x^n_{\lambda(v_jv_q)})=c_i$, a~contradiction. Since~$\mathcal{D}$ is acyclic, we obtain a~directed path~$P_i$ from~$s_i$ to~$t_i$.
    
        By the preceding argument, for each incoming edge $v_pv_j$, a~color~$c_i$ in position $\lambda(v_j)$ cannot be exchanged with the color in position $\lambda(v_pv_j)$. For each outgoing edge $v_jv_q$, position $\lambda(v_jv_q)$ still has its initial color $\sigma(x_{\lambda(v_jv_q)})\notin\{c_1,\ldots,c_\ell\}$ when the edge is considered at~$v_j$. Thus, exchanging the colors in positions $\lambda(v_j)$ and $\lambda(v_jv_q)$ gives a~color outside $\{c_1,\ldots,c_\ell\}$ in position $\lambda(v_j)$.
    
        Hence, if~$v_j\in V(P_i)$, then $\sigma(x^j_{\lambda(v_j)})=c_i$. Since~$c_1,\ldots,c_\ell$ are distinct, each vertex~$v_j$ belongs to at most~$1$ of the paths~$P_1,\ldots,P_\ell$. Therefore, the paths are vertex-disjoint.
    
        For each~$j\in[n]$, the clique~$A_j$ is obtained after the incoming operations at~$v_j$ and before its outgoing operations. By the preceding argument,~$v_j\in V(P_i)$ if and only if $\sigma(x^j_{\lambda(v_j)})=c_i$. Therefore, $V(P_i)=\{v_j\in V(\mathcal{D}):\sigma(x^j_{\lambda(v_j)})=c_i\}$ for each~$i\in[\ell]$.
    \end{clproof}
    
    We now show that~$G_{\mathcal{D}}$ is~$k$-colorable if and only if there exist pairwise vertex-disjoint directed paths~$P_1,\ldots,P_\ell$ in~$\mathcal{D}$ such that~$P_i$ is a~path from~$s_i$ to~$t_i$ for each~$i\in[\ell]$. We first assume that there exist pairwise vertex-disjoint directed paths from~$s_i$ to~$t_i$ for each~$i \in [\ell]$. By Claim~\ref{claim:DAGpathpropagate},~$H$ has a~proper~$k$-coloring~$\sigma$ such that $\sigma(x^n_{2i})=\sigma(x_{2i-1})$ for each~$i\in[\ell]$.
    
    Recall that~$R=[k-1]\setminus\{2,4,\ldots,2\ell\}$. In the construction step (2), all pairs~$p,q\in R$ with~$p<q$ are considered in lexicographic order, so the colors in these positions can be arranged in any order. We choose the exchanges so that position $2i-1$ has color $\sigma(x_{2i})$ for each~$i\in[\ell]$, and position~$r$ has color $\sigma(x_r)$ for each~$r\in\{2\ell+1,\ldots,k-1\}$. The color in position $2i$ remains $\sigma(x_{2i-1})$ for each~$i\in[\ell]$, and the color in position~$k$ remains $\sigma(x_k)$.
    
    In the construction step (3), each $\cs_{2i-1}$ exchanges the colors in positions $2i-1$ and $2i$. Thus, position~$r$ of~$A_{\mathrm{last}}$ has color $\sigma(x_r)$ for each~$r\in[k]$. By Lemma~\ref{lemma:coloringedgegadget},~$\sigma$ remains a~proper~$k$-coloring after adding the edges between~$A_{\mathrm{last}}$ and~$A_0$, and hence~$G_\mathcal{D}$ is~$k$-colorable.
    
    Conversely, suppose that~$G_\mathcal{D}$ has a~proper~$k$-coloring~$\sigma$, and write~$A_{\mathrm{last}}=\{y_1,\ldots,y_k\}$. By Lemma~\ref{lemma:coloringedgegadget}, the edges between~$A_{\mathrm{last}}$ and~$A_0$ give $\sigma(y_r)=\sigma(x_r)$ for each~$r\in[k]$. Since (2) keeps the colors in the sink positions unchanged and (3) exchanges positions $2i-1$ and $2i$, we obtain $\sigma(x^n_{2i})=\sigma(y_{2i-1})=\sigma(x_{2i-1})$ for each~$i\in[\ell]$. By Claim~\ref{claim:DAGpathpropagate}, there exist pairwise vertex-disjoint directed paths from~$s_i$ to~$t_i$ for each~$i\in[\ell]$.
    
    We next construct a~graph~$G_\mathcal{D}^*$ from~$G_\mathcal{D}$ such that~$G_\mathcal{D}^*$ is~$k$-colorable if and only if~$G_\mathcal{D}$ is~$k$-colorable. We construct~$G_\mathcal{D}^*$ from~$G_\mathcal{D}$ by inserting columns using operation $\ce$ and deleting the dummy vertex. Let~$B_0,\ldots,B_{N-1}$ be the columns of~$G_\mathcal{D}^*$ in cyclic order, where~$B_0=A_0\setminus \{x_k\}$. We choose the inserted columns as follows.
    \begin{enumerate}[(i)]
        \item For each original step~$B=f(A)$, where~$f\in \{\ce\}\cup \{\cc_j,\cs_j : j\in[k-1]\}$, we insert~$t \ge 1$ new columns~$C_1, \ldots, C_t$ between~$A$ and~$B$. Starting with~$C_0 = A$, we replace this step by successively obtaining~$C_r=\ce(C_{r-1})$ for each~$r\in[t]$, followed by~$B=f(C_t)$.
        
        \item We choose the number of $\ce$ operations before each original operation so that the following condition holds. 
        \begin{itemize}
            \item If~$B_r$ is introduced by $\cs_j$ or $\cc_j$, then~$r\equiv\nobreak j+1 \pmod{k-1}$.
        \end{itemize}
        This condition can be satisfied by inserting between $1$ and~$k-1$ $\ce$ operations.

        \item We insert at least~$1$ column between~$A_{\text{last}}$ and~$A_0$, using the edges of $\ce$ between consecutive columns. We choose at most~$k-1$ such columns so that~$N \equiv 1 \pmod{k-1}$.

        \item Finally, we delete the dummy vertex~$x_k$ in~$A_0$.
    \end{enumerate}

    Let~$B_0 = \{x_{0,1}, \ldots, x_{0, k-1}\}$ and let~$B_i = \{x_{i, 1}, \ldots, x_{i, k}\}$ for each~$i \in [N-1]$, where the vertices in each column are indexed according to their order in the construction.
    
    We show that~$G_\mathcal{D}$ is~$k$-colorable if and only if~$G_\mathcal{D}^*$ is~$k$-colorable. Suppose first that~$G_\mathcal{D}$ has a~proper~$k$-coloring. By Lemma~\ref{lemma:coloringedgegadget}, we obtain a~proper~$k$-coloring of~$G_\mathcal{D}^*$ by inserting the $\ce$ operations and then deleting~$x_k$.

    Conversely, suppose that~$G_\mathcal{D}^*$ has a~proper~$k$-coloring~$\sigma : V(G_\mathcal{D}^*) \to [k]$. Let~$c$ be the unique color in $[k]$ not used on~$B_0$. Consider~$B_1$ and~$B_{N-1}$. Since~$x_{1, k}$ is complete to~$B_0$, we obtain $\sigma(x_{1, k})=c$.

    Every $\cs_j$ or $\cc_j$ used in the construction satisfies~$j\le k-2$. Thus, by Lemma~\ref{lemma:coloringedgegadget}, the color in position~$k$ is unchanged along~$B_1,\ldots,B_{N-1}$, and hence $\sigma(x_{N-1, k})=c$. Since the neighbors of~$x_k$ before its deletion are precisely the vertices in~$B_0$ and~$B_{N-1}\setminus \{x_{N-1, k}\}$, none of these vertices has color~$c$. So, we can set $\sigma(x_k)=c$ to obtain a~proper~$k$-coloring of the graph before deletion. Along each inserted sequence of $\ce$ operations, corresponding positions have the same color. Thus, restricting~$\sigma$ to the original columns gives a~proper~$k$-coloring of~$G_\mathcal{D}$.

    Consequently,~$\mathcal{D}$ has the required disjoint paths if and only if~$G_\mathcal{D}^*$ is~$k$-colorable.

    We now show that $\tww(G_\mathcal{D}^*) \le 2$.

    \begin{claim}\label{claim:coloring-two-sequence}
        $\tww(G_\mathcal{D}^*) \le 2$.
    \end{claim}

    \begin{clproof}
        Let~$M = \frac{kN - 1}{k-1}$. We define the function~$C : V(G_\mathcal{D}^*) \to [kN-1]$ such that $C(x_{i, j}) = k(i + 1) - j$ for each~$x_{i, j} \in V(G_\mathcal{D}^*)$.

        To construct a~contraction sequence of~$G_\mathcal{D}^*$, we list the vertices in the order given by~$C$ and partition them into consecutive sets~$W_0,\ldots,W_{M-1}$ of size~$k-1$. For each~$p\in \{0,\ldots,M-1\}$, let~$W_p=\{w_{p,1},\ldots,w_{p,k-1}\}$, where~$w_{p,q}$ is the vertex~$x_{i,j}$ satisfying $C(x_{i,j})=(k-1)p+q$ for each~$q\in[k-1]$.

\begin{figure}[htbp]
\centering
\begin{tikzpicture}[font=\small,
  vertex/.style={circle,draw,fill=black,inner sep=1.5pt,outer sep=0pt},
  block/.style={draw,thick,rounded corners=4pt},
  regroup/.style={thick,dashed,rounded corners=3pt}]
  \draw[block,fill=blue!10] (-.36,-.48) rectangle (2.12,.50);
  \draw[block,fill=red!10] (2.28,-.48) rectangle (5.64,.50);
  \draw[block,fill=blue!10] (5.80,-.48) rectangle (9.16,.50);
  \draw[block,fill=red!10] (9.32,-.48) rectangle (12.68,.50);
  \foreach \x/\b in {.88/0,3.96/1,7.48/2,11/3}
    \node at (\x,.83) {$B_{\b}$};

  \foreach \c/\a/\b in {
    1/0/3,2/0/2,3/0/1,
    4/1/4,5/1/3,6/1/2,7/1/1,
    8/2/4,9/2/3,10/2/2,11/2/1,
    12/3/4,13/3/3,14/3/2,15/3/1}{
    \pgfmathsetmacro{\xx}{.88*(\c-1)}
    \node[vertex] at (\xx,0) {};
    \node at (\xx,.28) {\scriptsize $\c$};
    \node at (\xx,-.28) {$x_{\a,\b}$};
  }

  \foreach \p in {0,...,4}{
    \pgfmathsetmacro{\leftx}{2.64*\p-.30}
    \pgfmathsetmacro{\rightx}{2.64*\p+2.06}
    \pgfmathsetmacro{\midx}{2.64*\p+.88}
    \draw[regroup] (\leftx,-.59) -- (\leftx,-.91)
      -- (\rightx,-.91) -- (\rightx,-.59);
    \node at (\midx,-1.18) {$W_{\p}$};
  }
  \draw[->,thick] (-.32,1.24) -- (12.64,1.24);
  \node[fill=white,inner sep=3pt] at (6.16,1.24) {ordering~$C$};
\end{tikzpicture}
\caption{The vertices in the order given by~$C$ for~$k=N=4$.
Solid boxes show the columns~$B_i$, and dashed brackets show the cliques~$W_i$.
We omit all edges.}
\label{fig:coloring-reblocking}
\end{figure}
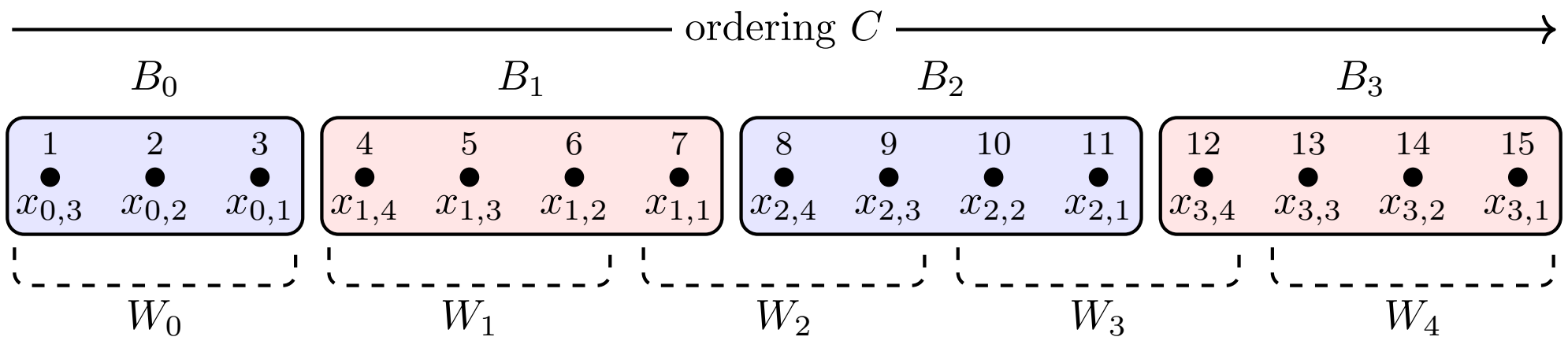

        First suppose that every operation is $\ce$. Apart from the closing edges, two vertices~$x_{a, b}$ and~$x_{a', b'}$ are adjacent if and only if $|C(x_{a, b}) - C(x_{a', b'})| \le k-1$. Thus, each~$W_i$ is a~clique. Let~$w_{i,p}=x_{r,a}$ and~$w_{i+1,q}=x_{s,b}$. By the definition of~$W_i$, we have $C(x_{s,b})-C(x_{r,a})=(k-1)+q-p$. Thus, if every operation is $\ce$, then~$w_{i, p}$ and~$w_{i+1, q}$ are adjacent if and only if~$q\le p$.

        Now suppose that~$B_r$ is introduced by $\cs_j$ or $\cc_j$. Both operations omit $x_{r-1,j}x_{r,j+1}$. Since~$r \equiv\nobreak j+1 \pmod{k-1}$, we have $C(x_{r-1,j})=kr-j=(k-1)i+1$ for some~$i$. Since $C(x_{r,j+1})=kr-j+k-1 = (k-1)(i+1) + 1$, the omitted edge is $w_{i,1}w_{i+1,1}$. For an operation $\cs_j$, the added edge is $x_{r-1,j}x_{r,j}$. Since $C(x_{r, j}) = kr - j +k = (k-1)(i+1) + 2$, this edge is $w_{i,1}w_{i+1,2}$.

        Consequently, each~$W_i$ forms a~clique, and the edges between consecutive columns are described as follows. Let~$F_i=\{w_{i,p}w_{i+1,q}:1\le q\le p\le k-1\}$. Then,
        \[
        E(W_i,W_{i+1}) =
        \begin{cases}
            F_i
            & \text{for } \ce, \\[4pt]
            \bigl(F_i \setminus \{w_{i,1}w_{i+1,1}\}\bigr)
            \cup \{w_{i,1}w_{i+1,2}\}
            & \text{for } \cs_j, \\[4pt]
            F_i \setminus \{w_{i,1}w_{i+1,1}\}
            & \text{for } \cc_j.
        \end{cases}
        \]
        
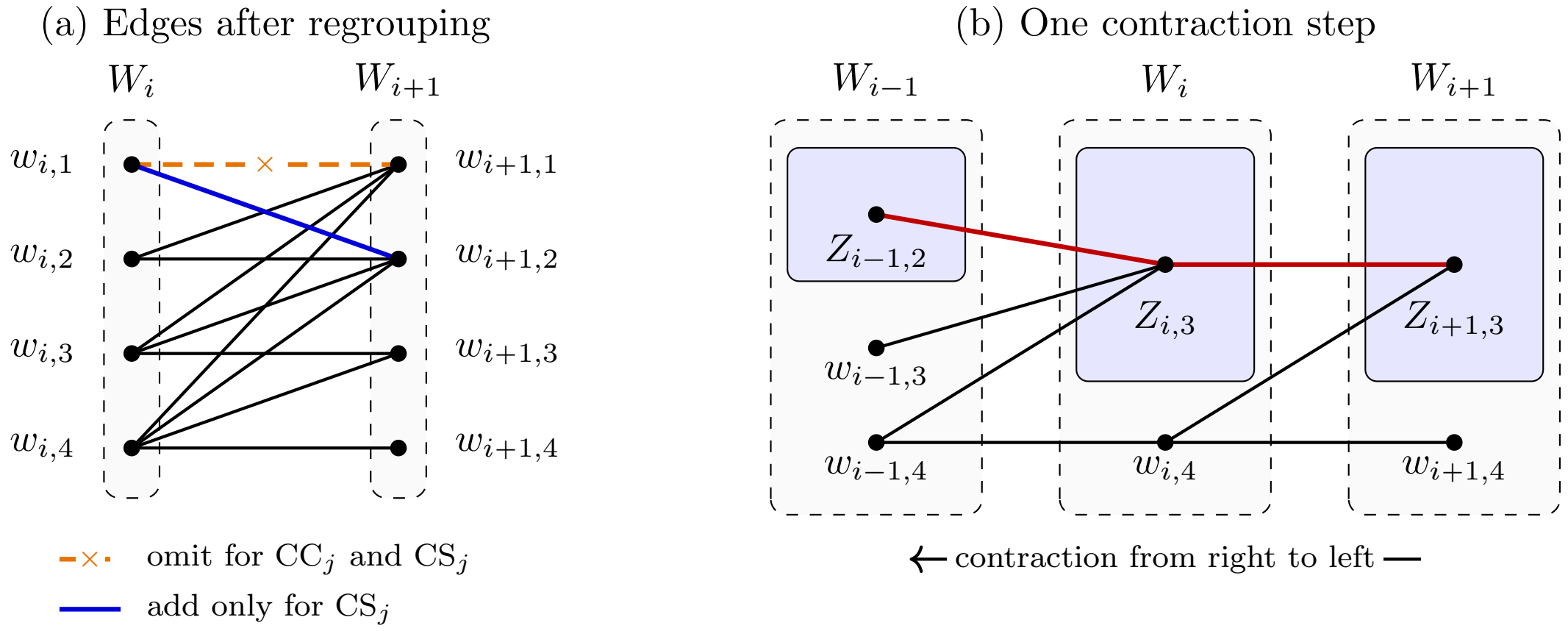
\begin{figure}[tbp]
\centering
\begin{tikzpicture}[font=\small,
  vertex/.style={circle,fill=black,inner sep=1.5pt},
  column/.style={draw,dashed,rounded corners=4pt,fill=gray!4},
  prefixpart/.style={draw,rounded corners=3pt,fill=blue!10},
  ordinary/.style={thick},
  missing/.style={orange!90!black,dash pattern=on 5pt off 3pt,line width=1.2pt},
  added/.style={blue!85!black,solid,line width=1.2pt},
  rededge/.style={red!75!black,solid,very thick}]
  \begin{scope}
    \node at (2.1,2.45) {(a) Edges after regrouping};
    \draw[column] (.65,-1.8) rectangle (1.15,1.6);
    \draw[column] (3.05,-1.8) rectangle (3.55,1.6);
    \node at (.9,1.95) {$W_i$};
    \node at (3.3,1.95) {$W_{i+1}$};
    \foreach \p in {1,...,4} {
      \coordinate (L\p) at (.9,{2.05-.85*\p});
      \coordinate (R\p) at (3.3,{2.05-.85*\p});
    }
    \foreach \p in {1,...,4} {
      \foreach \q in {1,...,\p} {
        \ifnum\p=1\relax
        \else
          \draw[ordinary] (L\p)--(R\q);
        \fi
      }
    }
    \draw[missing] (L1)--(R1);
    \node[text=orange!90!black,fill=white,inner sep=.4pt,font=\scriptsize]
      at (2.1,1.2) {$\times$};
    \draw[added] (L1)--(R2);
    \foreach \p in {1,...,4} {
      \node[vertex] at (L\p) {};
      \node[vertex] at (R\p) {};
      \node[anchor=east,inner sep=0pt] at (.40,{2.05-.85*\p}) {$w_{i,\p}$};
      \node[anchor=west,inner sep=0pt] at (3.80,{2.05-.85*\p}) {$w_{i+1,\p}$};
    }
    \draw[missing] (.25,-2.35)--(.8,-2.35);
    \node[text=orange!90!black,fill=white,inner sep=.4pt,font=\scriptsize]
      at (.525,-2.35) {$\times$};
    \node[anchor=west,font=\scriptsize] at (.9,-2.35) {omit for $\cc_j$ and $\cs_j$};
    \draw[added] (.25,-2.8)--(.8,-2.8);
    \node[anchor=west,font=\scriptsize] at (.9,-2.8) {add only for $\cs_j$};
  \end{scope}
  \begin{scope}[xshift=6.7cm]
    \node at (3.5,2.45) {(b) One contraction step};
    \foreach \x/\lab in {.9/{W_{i-1}},3.5/{W_i},6.1/{W_{i+1}}} {
      \draw[column] ({\x-.95},-1.95) rectangle ({\x+.95},1.6);
      \node at (\x,1.95) {$\lab$};
    }
    \draw[prefixpart] (.1,.15) rectangle (1.7,1.35);
    \draw[prefixpart] (2.7,-.75) rectangle (4.3,1.35);
    \draw[prefixpart] (5.3,-.75) rectangle (6.9,1.35);
    \coordinate (PL) at (.9,.75);
    \coordinate (PC) at (3.5,.30);
    \coordinate (PR) at (6.1,.30);
    \coordinate (L3) at (.9,-.45);
    \coordinate (L4) at (.9,-1.3);
    \coordinate (C4) at (3.5,-1.3);
    \coordinate (R4) at (6.1,-1.3);
    \draw[ordinary] (L3)--(PC);
    \draw[ordinary] (L4)--(PC);
    \draw[ordinary] (L4)--(C4);
    \draw[ordinary] (C4)--(PR);
    \draw[ordinary] (C4)--(R4);
    \draw[rededge] (PL)--(PC)--(PR);
    \node at (.9,.40) {$Z_{i-1,2}$};
    \node at (3.5,-.20) {$Z_{i,3}$};
    \node at (6.1,-.20) {$Z_{i+1,3}$};
    \foreach \v in {PL,PC,PR} {\node[vertex] at (\v) {};}
    \foreach \v in {L3,L4,C4,R4} {\node[vertex] at (\v) {};}
    \node[anchor=north,inner sep=0pt] at (.9,-.57) {$w_{i-1,3}$};
    \node[anchor=north,inner sep=0pt] at (.9,-1.42) {$w_{i-1,4}$};
    \node[anchor=north,inner sep=0pt] at (3.5,-1.42) {$w_{i,4}$};
    \node[anchor=north,inner sep=0pt] at (6.1,-1.42) {$w_{i+1,4}$};
    \draw[->,thick] (5.8,-2.35)--(1.2,-2.35);
    \node[font=\scriptsize,fill=white,inner sep=2pt] at (3.5,-2.35)
      {contraction from right to left};
  \end{scope}
\end{tikzpicture}
\caption{Consecutive columns for~$k=5$ in (a) and one contraction with~$j=3$ in (b).
Orange and blue lines show omitted and added edges, respectively.
We omit edges within each column.}
\label{fig:coloring-sequence}
\end{figure}

        Finally,~$w_{M-1,p}=x_{N-1,k-p}$ and~$w_{0,q}=x_{0,k-q}$. The closing $\ce$ condition is~$k-p<k-q$, so $E(W_{M-1},W_0)=\{w_{M-1,p}w_{0,q}:1\le q<p\le k-1\}$. There are no edges between other pairs of distinct columns.

        Now, we construct a~contraction sequence of~$G_\mathcal{D}^*$. For each~$i\in \{0,\ldots,M-1\}$ and~$j\in[k-1]$, let~$Z_{i,j}=\{w_{i,1},\ldots,w_{i,j}\}$. For~$j=2,\ldots,k-1$, we consider~$i=M-1,\ldots,0$ in this order and contract~$Z_{i,j-1}$ and $\{w_{i,j}\}$ to obtain~$Z_{i,j}$.

        We show that every part has red degree at most~$2$. Consider the contraction obtaining~$Z_{i,j}$ with~$j \ge 2$. If $0<i<M-1$, then~$Z_{i-1,j-1}$ and~$Z_{i+1,j}$ are parts. By the construction, every vertex~$w_{i-1,p}$ with~$p\ge j$ is complete to~$Z_{i,j}$, and every vertex~$w_{i+1,q}$ with~$q>j$ is anti-complete to~$Z_{i,j}$. Thus, the only possible red neighbors of~$Z_{i,j}$ are~$Z_{i-1,j-1}$ and~$Z_{i+1,j}$.

        If~$i=M-1$, then~$Z_{0,j-1}$ is a~part. Since every closing edge $w_{M-1, a}w_{0, b}$ satisfies~$b < a$, every vertex~$w_{0, q}$ with~$q \ge j$ is anti-complete to~$Z_{M-1, j}$. Thus, its only possible red neighbors are~$Z_{M-2, j-1}$ and~$Z_{0, j-1}$.

        If~$i=0$, then~$Z_{M-1,j}$ has already been obtained. Every vertex~$w_{M-1,p}$ with~$p>j$ is complete to~$Z_{0,j}$. Thus, its only possible red neighbors are~$Z_{M-1,j}$ and~$Z_{1,j}$.

        By the construction, there are no edges between nonconsecutive columns. Therefore, every part has red degree at most~$2$.

        Finally, the parts~$Z_{0,k-1},\ldots,Z_{M-1,k-1}$ form a~red cycle. We successively contract adjacent parts of this cycle until only one part remains, keeping the red degree at most~$2$.
    \end{clproof}

    Each $\cc^*_{(p,q)}$ introduces $O(k)$ columns, and we use $O(k^2)$ such operations. Hence,~$G_\mathcal{D}$ has $O(k^3)$ columns. We insert at most~$k-1$ columns before each original operation and at most~$k-1$ columns between~$A_{\text{last}}$ and~$A_0$, so~$N=O(k^4)$ and $|V(G_\mathcal{D}^*)|=O(k^5)$.

    As shown above, the required pairwise vertex-disjoint directed paths exist in~$\mathcal{D} = (V, E)$ if and only if~$G_\mathcal{D}^*$ is $(|V| + |E| + 1)$-colorable. Moreover, the above contraction sequence witnesses that $\tww(G_\mathcal{D}^*)\le 2$, and both~$G_\mathcal{D}^*$ and the contraction sequence can be constructed in polynomial time. Thus, this gives a~polynomial-time reduction from \textsc{Disjoint Paths on Directed Acyclic Graphs} to \textsc{Minimum Coloring} on graphs of twin-width at most~$2$, even when a~witnessing contraction sequence is given. Since a~proper~$k$-coloring can be verified in polynomial time, \textsc{Minimum Coloring} is NP-complete on this class.
\end{proof}

\paragraph*{AI disclosure.}
The authors first obtained a~mim-width bound of $105$ for graphs of twin-width at most~$2$, and improved it to $8$ with the assistance of GPT-6.1 Sol.
They also used GPT-6.1 Sol to improve their linear bound on the tree-independence number.

The structural idea for the \textsc{Minimum Coloring} reduction is due to the authors.
GPT-6.1 Sol assisted with the construction ensuring that the resulting graphs have twin-width at most~$2$.

GPT-6.1 Sol also helped review the manuscript, correct typos, adjust the page layout, and add figures.

\bibliographystyle{plain}
\bibliography{mim}

\end{document}